\documentclass[11pt]{article}

\usepackage{latexsym}

\font\tenbb=msbm10
\font\sevenbb=msbm7
\font\fivebb=msbm5
\newfam\bbfam
\textfont\bbfam=\tenbb
\scriptfont\bbfam=\sevenbb
\scriptscriptfont\bbfam=\fivebb
\def\Bbb{\fam\bbfam\tenbb}

\newcommand{\begeq}[1]{\begin{equation} \label{#1}}
\newcommand{\fineq}{\end{equation}}

\newcommand{\ra}{\rightarrow}

\newcommand{\m}{\displaystyle}

\newcommand{\imp}{\Rightarrow}
\newcommand{\g}{\gamma}
\newcommand{\al}{\alpha}
\newcommand{\e}{\epsilon}
\newcommand{\de}{\delta}
\newcommand{\ka}{\kappa}
\newcommand{\ti}{\tilde}

\def\bR{{\Bbb R}}                  
\def\bE{{\Bbb E}}                  
\def\bP{{\Bbb P}}                  

\def\tinf{{\rightarrow\infty}}     

\def\tts{{\textstyle}}

\begin{document}

\newtheorem{theo}{Theorem}
\newtheorem{prop}{Proposition}
\newtheorem{lemma}{Lemma}
\newtheorem{cor}{Corollary}
\newtheorem{defi}{Definition}

\begin{titlepage}

\title{Confidence bands for densities, logarithmic point of view}
\author{Abdelkader Mokkadem \and Mariane Pelletier
}
\date{\small{
Universit\'e de Versailles-Saint-Quentin\\
D\'epartement de Math\'ematiques\\
45, Avenue des Etats-Unis\\
78035 Versailles Cedex\\
France\\}
(mokkadem,pelletier)@math.uvsq.fr\\}
\maketitle

\begin{abstract}
Let $f$ be a probability density and $C$ be an interval on which $f$ is
bounded away from zero. By establishing the limiting distribution of the 
uniform error of the kernel estimates $f_n$ of $f$, 
Bickel and Rosenblatt (1973) 
provide confidence bands $B_n$ for $f$ on $C$ with asymptotic level
$1-\al\in]0,1[$. Each of the confidence intervals whose 
union gives
$B_n$ has an asymptotic level equal to one; pointwise moderate 
deviations principles 
allow to prove that all these intervals share the same logarithmic
asymptotic level. Now, as soon as both pointwise and uniform moderate 
deviations principles for $f_n$ exist, they share the same asymptotics. 
Taking this observation as a starting point, we present a new approach 
for the construction of confidence bands for $f$, based on the use 
of moderate deviations principles. The advantages of this approach are 
the following: (i) it enables to construct confidence bands, which have 
the same width (or even a smaller width) as the confidence bands 
provided by Bickel and Rosenblatt (1973), but which have a better 
aymptotic level; (ii) any confidence band constructed in that way 
shares the same logarithmic asymptotic level as all the confidence 
intervals, which make up this confidence band;  
(iii) it allows to deal with all the dimensions in 
the same way; (iv) it enables to sort out the problem 
of providing confidence bands for $f$ on compact sets on 
which $f$ vanishes (or on all $\bR^d$), by introducing 
a truncating operation. 
\end{abstract}
\vspace{2cm}


\end{titlepage}

\section{Introduction}

Let $X_1,\ldots,X_n$ be independent and identically distributed 
$\bR^d$-valued random variables with bounded probability density 
function $f$. The problem of computing confidence 
bands for $f$ is of central interest in 
nonparametric statistics. \\

One main known approach to the construction of
confidence bands for $f$ is due to Bickel and 
Rosenblatt (1973) in the case $d=1$ and to Rosenblatt (1976) 
in the case $d\geq 2$, and is based on the 
limiting distribution of the normalized uniform error of the kernel 
density estimate.
The approach of Bickel and Rosenblatt has been 
extended, in the unidimensional case, in several directions; 
among others, let us cite Mack (1982) and Liu and Ryzin (1986) 
for the extension to other types of density estimates, 
Burke and Horvath (1984) and Mielniczuk (1987) for the 
censored data case, Xu and Martinsek (1995), Martinsek and
Xu (1996), Sun and Zhou (1998) for the construction of sequential 
confidence bands, and Gin\'e, Koltchinskii and Sakhanenko (2003, 2004) 
for different asymptotics of the weighted uniform error.  
Another common technique for constructing confidence bands is 
through the bootstrap, which is in particular used for bias estimation; 
see, for example, Hall (1992) for the density 
and H\"{a}rdle and Marron (1991) for the regression. 
For other approaches, see Hall and Titterington (1988) and
Hall and Owen (1993).\\

Our object in this paper is to present a new approach, based on the use 
of Moderate Deviations Principles (MDP) of the normalized uniform 
error of the kernel density estimates. We avoid bias estimation by a 
slight undersmoothing, which is shown in 
Hall (1992) to be more efficient than explicit bias 
correction when the goal is to minimize the coverage error of the 
confidence band. 

The use of large and moderate deviations in statistical inference 
is not new. It has been initiated by the papers of 
Chernoff (1952) and Bahadur (1960), and then 
developped in various directions. Let us cite, among many others,  
Borovkov and Mogulski (1992), Groeneboom (1980), 
Ibragimov and Radavicius (1981), Kallenberg (1982, 1983a, 1983b), 
Korostelev and Leonov (1995), Nikitin (1995), Mokkadem and Pelletier 
(2005), and Puhalskii and Spokoiny (1998).\\

The idea of using MDP
for the construction of confidence bands for $f$ comes naturally 
when making a pointwise analysis of the confidence bands provided  
by Bickel and Rosenblatt (1973).  
Consider the univariate framework (that is, the case $d=1$), and let 
$C$ be a bounded interval of $\bR$ on which $f$ is assumed 
to be bounded away from zero. Let $f_n$ denote the kernel estimator of
$f$;  
Bickel and Rosenblatt (1973) establish the asymptotic law of 
$\sup_{x\in C}\left|f_n(x)-f(x)\right|/\sqrt{f(x)}$ 
suitably normalized. 
This result allows them to provide sequences of random intervals 
$I_n(x)$ for all $x\in C$, which satisfy the property:
$$\lim_{n\tinf}\bP\left(\exists x\in C,\ f(x)\not\in I_n(x)\right)=\al.$$
For simplicity, we also denote by $I_n(x)$ the segment 
$\{x\}\times I_n(x)$, and say that 
$$B_{n,\al}=\cup_{x\in C}I_n(x)$$
is a confidence band for $f$ on $C$ with asymptotic level 
$1-\al\in]0,1[$. 
In other words, the set of functions
$$D_{n,\al}=\{g:\bR\rightarrow\bR ,\ g(x)\in I_n(x)\ \forall x\in C\}$$ 
is a confidence region of $f$ with asymptotic level $1-\al$. (Although 
$D_{n,\al}$ is a confidence region of the functional parameter $f$ since
$\lim_{n\tinf}\bP\left(f\in D_{n,\al}\right)=1-\al$, it gives nontrivial 
confidence intervals of the values $f(x)$ only for $x\in C$).  

Now, a straightforward application of the central limit theorem (CLT) allows 
to prove that the asymptotic level of each confidence interval $I_n(x)$ 
($x\in C$) whose union gives $B_{n,\alpha}$, is one. A natural question 
is then to wonder at what rate the levels of the intervals $I_n(x)$ 
go to one. A result giving the convergence rate 
to zero of the sequence 
$\bP\left(f(x)\not\in I_n(x)\right)$ is typically a  
MDP result; this convergence rate is thus 
expected to be exponential. That is the reason why we introduce 
here the notion of 
\emph{logarithmic asymptotic level} for confidence regions of (eventually
infinite dimensional) unknown parameters.

\begin{defi}
Let $\{D_n\}$ be a sequence of confidence regions of an unknown parameter
$\theta$. The logarithmic asymptotic level of
 $\{D_n\}$ is $\gamma$ ($\gamma>0$) with speed $w_n$ ($w_n\ra\infty$) if
$$\lim_{n\tinf}\frac{1}{w_n}\log\bP(\theta\not\in D_n)=-\gamma.$$
\end{defi}
Of course, if $\{D_n\}$ has a logarithmic asymptotic level $\gamma>0$, then
the asymptotic level of $\{D_n\}$ is necessarily one.\\

It turns out that the sequences of confidence regions, which have 
a positive
logarithmic asymptotic level, are often asymptotic almost sure sequences
of confidence regions in the sense of the following definition.

\begin{defi}
Let $\{D_n\}$ be a sequence of confidence regions of an unknown parameter
$\theta$, and let $\Omega$ denote the underlying probability space.
$\{D_n\}$ is an asymptotic almost sure (or consistent) 
sequence of confidence regions 
of $\theta$ if there exists $\Omega_0\subset\Omega$ such that: 
$$\left\{\begin{array}{ll}
\bullet & \bP(\Omega_0)=1 \\
\bullet & \forall\omega\in\Omega_0, \ \exists N(\omega)\  \mbox{such that}\ 
n\geq N(\omega)\imp\theta\in D_n(\omega).
\end{array}\right.
$$
\end{defi}
\vspace{0.5cm}

Indeed, the following proposition is a straightforward consequence of 
Borel and Cantelli Lemma. 

\begin{prop}
\label{aspro}
Let $\{D_n\}$ be a sequence of confidence regions of an unknown parameter
$\theta$, whose logarithmic asymptotic level is $\gamma>0$ with speed $w_n\ra \infty$. 
If there exists $\delta\in]0,\g[$ such that $\sum \exp(-\delta w_n)<\infty$, then 
$\{D_n\}$ is an asymptotic almost sure sequence of confidence regions of 
$\theta$.
\end{prop}

Let us come back to the pointwise analysis of the confidence bands 
provided by Bickel and Rosenblatt (1973). Later on, we shall prove 
in particular that, 
for all $x\in C$, the intervals $I_n(x)$ have a logarithmic
asymptotic level 1 with speed ${\log(1/h_n)}$, where $h_n$ is the 
bandwidth used for the computation of the kernel estimator $f_n$. 
So, the confidence
bands $B_{n,\al}$ with asymptotic level $1-\al<1$ provided by
Bickel and Rosenblatt are unions of confidence intervals
$I_n(x)$ whose asymptotic levels equal one and whose logarithmic 
asymptotic levels are independent on the value of $x\in C$. \\

The difference between the asymptotic level of the confidence band 
$B_{n,\al}$ on the one hand and the asymptotic 
levels of all the confidence 
intervals $I_n(x)$ on the other hand is explained by the difference 
between the asymptotic weak behaviour of the uniform error 
of the kernel density estimator (given by Bickel and 
Rosenblatt's result) on the one hand 
and the asymptotic weak behaviour of the pointwise error of the kernel 
density estimator (given by the central limit theorem) on 
the other hand. 
Now, MDP for the nonnormalized error of the kernel density 
estimator have been established by Gao (2003) (see also Mokkadem, 
Pelletier and Worms (2005)); it turns out that, as soon as 
both pointwise and uniform MDP exist, the pointwise
and the uniform MDP share exactly the same asymptotics. Taking this 
remark as a starting point, we propose, in this paper, a new approach 
to construct confidence bands for $f$ based on the use of MDP for 
the normalized error of the kernel density estimator.
This approach has several advantages:
\begin{itemize}
\item It allows to construct confidence bands $B^*_n$, which have 
the same width (or even a smaller width) as the confidence bands 
$B_{n,\al}$ provided by Bickel and Rosenblatt (1973), but which:
(i) have an asymptotic level equal to one instead of 
$1-\al\in]0,1[$; (ii) share the same logarithmic asymptotic level as 
all the confidence intervals whose union gives $B^*_n$; 
(iii) are asymptotic almost sure confidence bands.
\item In order to deal with the multivariate framework, 
Rosenblatt (1976) has to require the use of higher order kernels and, 
consequently, to impose rather stringent conditions on $f$; 
in contrast, in the MDP approach, all the dimensions are dealt with in 
the same way, and thus without any additional 
assumption neither on the density, nor on the kernel, 
in the case $d\geq 2$. 
\item Whatever the dimension $d$ is, Bickel and Rosenblatt require the 
condition ``$f$ is bounded away from zero on $C$''. On the 
contrary, our approach enables us to sort out the problem 
of providing confidence bands for $f$ on compact sets on 
which $f$ vanishes. As a matter of fact, we introduce 
a truncating operation, which modifies the width of our 
confidence bands at some points $x\in C$, but which does not 
affect the logarithmic asymptotic level of our 
confidence bands. This truncating operation also enables 
us to provide confidence bands for $f$ on all $\bR^d$.
Let us mention that, in the case $d=1$, Gin\'e, Koltchinskii and 
Sakhanenko (2003, 2004) 
propose a slight modification of  Bickel and Rosenblatt's 
normalization of the uniform error; this allows them to construct confidence bands 
on the whole line in the case $f$ does not vanish on $\bR$.
\end{itemize}

Our paper is now organized as follows. In Section 2, we explicit 
the construction of our confidence bands. Section 3 is devoted 
to the precise statement of our assumptions and main results.
In Section 4, we discuss particular examples of applications of 
our main results: we first come back on the pointwise analysis 
of the confidence bands provided by Bickel and Rosenblatt (1973), 
and show how the MDP approach allows to construct more suitable 
confidence bands; then, we consider the problem of constructing confidence bands with smaller width. 
Section 5 is reserved to the proofs.

\section{Construction of confidence bands based on the use of MDP}

Let $C$ be a subset of $\bR^d$ and $(v_n)$ be a positive nonrandom 
sequence that goes to infinity. In this section, we construct 
confidence bands for $f$ on $C$ with width of order $v_n^{-1}$. 
We first consider the case $C$ is a compact set on which $f$ is 
bounded away from zero, and then introduce a truncating operation, 
which allows to consider the general framework ($f$ may vanish on 
$C$, $C$ may equal $\bR^d$).

\subsection{Confidence bands on compact sets on which $f$ is 
bounded away from zero}

Let $C$ be a compact set of $\bR^d$ on which $f$ is bounded 
away from zero. 
To construct a confidence band for $f$ on $C$ with width of order 
$v_n^{-1}$, we first construct, for all $x\in C$, confidence 
intervals of $f(x)$ with width of the same order. For that 
purpose, we proceed as follows.
\begin{itemize}
\item 
We estimate $f(x)$ by using the kernel estimator 
\begin{equation}
\label{g1}
f^*_n(x) = \frac{1}{nh_n^{*d}} \sum_{j=1}^n K\left( 
\frac{x-X_j}{h_n^*}\right),
\end{equation}
where the bandwidth $(h_n^*)$ is a sequence of positive 
real numbers such that
$\lim_{n\tinf} h_n^{*} =0$, $\lim_{n\tinf} nh_n^{*d} = \infty$, 
and where the kernel $K$ is a bounded nonnegative function
satisfying $\int_{\bR^d} K(z)\, dz =1$ and 
$\lim_{\|z\|\tinf} K(z) = 0$. 
\item 
The variance of $f^*_n(x)$ is equivalent (as $n$ goes to infinity) 
to $(nh_n^{*d})^{-1}f(x)\kappa$ where $\kappa=\int_{\bR^d} K^2(x)dx$; 
we estimate it by 
$(nh_n^{*d})^{-1}f_n(x)\kappa$, where $f_n$ is the kernel estimator 
of $f$ defined by 
\begin{equation}
\label{g2}
 f_n(x) = \frac{1}{nh_n^d} \sum_{j=1}^n K\left( \frac{x-X_j}{h_n}\right),
\end{equation}
the bandwidth $(h_n)$ being eventually different from $(h_n^*)$.
\item 
The confidence intervals for $f(x)$ ($x\in C$) are then defined as 
\begin{equation}
\label{Ihat}
\hat I_{n}(x)=  
\left[f_n^*(x)-\delta\frac{\sqrt{f_n(x)\kappa}}{v_n}\ ;\
f_n^*(x)+\delta\frac{\sqrt{f_n(x)\kappa}}{v_n}\right],
\end{equation}
where $\de>0$.
\end{itemize}
Our confidence band for $f$ on $C$ is finally defined by setting:
\begin{equation}
\label{Bhat}
\hat B_n=\cup_{x\in C}\hat I_n(x).
\end{equation}
In Section \ref{2AMA}, we give conditions on the sequence $(v_n)$ 
and the bandwidths $(h_n)$ and $(h_n^*)$, which ensure 
that the logarithmic asymptotic level of 
each interval $\hat I_n(x)$, $x\in C$, on the one hand,  
and of the confidence band $\hat B_n$ on the other hand,  
is $\de^2/2$ with speed $nh_n^{*d}/v^2_n$ (see Theorems 
\ref{th0} and \ref{th1}).

\subsection{Truncating operation} 
 
In order to allow the construction of confidence bands for 
$f$ on subsets $C$ of $\bR^d$ (eventually equal to $\bR^d$) 
on which $f$ may take the value zero, we now introduce a
truncating method, which relies on the following fact.  
For the values of $x\in C$ for which $f_n(x)$ is ``large 
enough'', the width of the intervals $\hat I_n(x)$ defined 
in (\ref{Ihat}) is suitable; 
but, for the values of $x\in C$ for which $f_n(x)$ is
zero, or, more generally, ``close to zero'', the width of the 
intervals $\hat I_n(x)$ is clearly not appropriate any more. 
In order to compensate for this problem which appears for 
``small'' values of $f_n(x)$, we impose a minimum width 
to all the confidence intervals whose union gives the confidence
band for $f$ on $C$; of course, this minimum width does not 
affect the width of the confidence band at the points 
$x\in C$ for which $f_n(x)$ is ``large enough''.\\

More precisely, we introduce a sequence $(\e_n)$ of positive 
real numbers satisfying $\lim_{n\tinf}\e_n= 0$, and define the 
truncating function $\tilde T_n$ by setting
\begin{equation}
\label{tildtdef}
\tilde T_n(x)=\left\{\begin{array}{ll}
f_n(x) & \mbox{if } f_n(x)\geq \epsilon_n, \\
\epsilon_n & \mbox{otherwise.} 
\end{array}\right.
\end{equation}
For each $x\in C$, we set 
\begin{equation}
\label{Icheck}
\check{I}_n(x)=\left[f_n^*(x)-\de
\frac{\sqrt{{\tilde T}_n(x)\kappa}}{v_n}\ ;\
f_n^*(x)+\de
\frac{\sqrt{{\tilde T}_n(x)\kappa}}{{v_n}}\right]
\end{equation}
and finally define $\check{B}_n$ as
\begin{equation}
\label{Bcheck}
\check{B}_n=\cup_{x\in C}\check{I}_n(x).
\end{equation}
In Section \ref{3AMA}, we give conditions on $(\e_n)$, which ensure 
that the logarithmic asymptotic level of the confidence band 
$\check{B}_n$ is $\de^2/2$ with speed $nh_n^{*d}/v^2_n$ (see 
Theorem \ref{th2} in the case $C$ is a compact set, and 
Theorem \ref{th3} in the case $C=\bR^d$).
In other words, the logarithmic asymptotic level of the 
confidence band $\hat{B}_n$ defined in (\ref{Bhat})
is not affected by the introduction of this truncating method.\\

From a practical point of view, it seems more realistic to 
take the width of the largest confidence interval into account
in the truncating operation, that is, to introduce the quantity 
$\sup_{x\in C}f_n(x)$ in the definition of the truncating function. 
For that purpose, let the sequence $(\e_n)$ satisfy the additional 
condition $\e_n\leq 1$ for all $n$, and define 
the function $T_n$ by setting
\begin{equation}
\label{Tdef}
T_n(x)=\left\{\begin{array}{ll}
f_n(x) & \mbox{if } f_n(x)\geq \epsilon_n [\sup_{x\in C}f_n(x)]\\
\epsilon_n[\sup_{x\in C}f_n(x)] & \mbox{otherwise.} 
\end{array}\right.
\end{equation}
In Section \ref{3AMA}, we establish that when the parameter 
$\tilde T_n(x)$ is replaced by $T_n(x)$ in the intervals 
$\check{I}_n(x)$, the logarithmic asymptotic level of the 
confidence band $\check{B}_n$ remains unchanged (see 
Corollary \ref{co1} in the case $C$ is a compact set, and 
Corollary \ref{co2} in the case $C=\bR^d$).

\section{Assumptions and Main Results}
\label{AMA}

\subsection{Assumptions}
\label{1AMA}
Before stating our assumptions, let us first define the 
\emph{covering number condition}.
Let $Q$ be a probability on $\bR^d$ and ${\mathcal F}\subset{\cal L}_2(Q)$ 
be a class of
$Q$-integrable functions. The covering number (see
Pollard (1984)) is the smallest value  $N_2(\e ,Q,{\cal F})$
of $m$ for which there exist $m$ functions $g_1,\ldots ,g_m\in {\cal L}_2(Q)$ 
such that
$$\min_{ i\in\{1,\ldots ,m\}}\left\|f-g_i\right\|_{{\cal L}_2(Q)}
\leq\e\ \ \forall f\in{\cal F}$$
(if no such $m$ exists, $N_2(\e ,Q,{\cal F})=\infty$).
Now, let $\Lambda$ be a bounded and integrable function on $\bR^d$, and 
let ${\cal F}(\Lambda)$ be the class of functions defined by
\begin{equation}
\label{Fcovering}
{\cal F}(\Lambda)=\left\{z\mapsto
\Lambda\left(\frac{x-z}{h}\right),\ \ h>0,\ \ x\in\bR^d\right\}.
\end{equation}
\emph{$\Lambda$ is said to satisfy the covering number condition if there exist
$A>0$ and $v>0$ such that, for any probability $Q$
on $\bR^d$ and any $\e\in]0,1[$,
\begin{equation}
\label{covering}
N_2(\e\|\Lambda\|_\infty ,Q,{\cal F}(\Lambda))\leq 
{\left(\frac{A}{\e}\right)}^{v}.
\end{equation}
}
\\
The classes which satisfy (\ref{covering}) are often called 
Vapnik-Chervonenkis classes. When $d=1$, the
real valued kernels with bounded variations satisfy the covering number condition (see
Pollard (1984)). Some examples of multivariate kernels satisfying the covering number
condition are the following~:

- the kernels defined as $K(x)=\psi(\|x\|)$, where $\psi$ is a real valued function
with bounded variations (see Nolan and Pollard (1987)).

- the kernels defined as $K(x)=\prod_{i=1}^dK_i\left(x_i\right)$ where the $K_i$,
$1\leq i\leq d$, are real valued functions with bounded variations (this follows from
Lemma A1 in Einmahl and Mason (2000)).\\
\\

We can now state our assumptions.\\

$X_1,\ldots ,X_n$ are \emph{i.i.d.} $\bR^d$-valued random vectors with 
bounded probability density $f$.
The kernel estimators $f_n$ and $f^*_n$ of $f$ are defined in (\ref{g2}) and 
(\ref{g1}) respectively, and the bandwidths 
$(h_n)$ and $(h_n^*)$ are two sequences of positive real numbers such 
that
$$h_n\rightarrow 0\ \ \mbox{and}\ \ h_n^*\rightarrow 0.$$

The assumptions to which we will refer in the sequel are the following.
\begin{itemize}
\item[$(A1)$] $K$ is a bounded and nonnegative function on $\bR^d$
such that
\[
 \textstyle \int_{\bR^d} K(z)dz =1, \ \ \int_{\bR^d} z_jK(z)dz=0
\ \forall j\in\{1,\dots ,d\},\ \ \mbox{and} 
\ \ \int_{\bR^d} \|z\|^2|K(z)|dz < \infty.
\]  

\item[$(A2)$]  $K$ is H\"older-continuous on $\bR^d$ and satisfies the 
covering number condition.

\item[$(A3)$] $f$ is twice differentiable on $\bR^d$,
$\sup_{x\in\bR^d}\|\nabla f(x)\|<\infty$, and 
$\sup_{x\in\bR^d}\|D^2f(x)\|<\infty$.   
\item[$(A4)$]
There exists $q>0$ such that $z\mapsto \|z\|^{q} f(z)$ is a 
bounded function on $\bR^d$.
\end{itemize}
\medskip
Let us recall the notation
$$\kappa=\int_{\bR^d}K^2(z)dz.$$

\subsection{Confidence regions without truncating}
\label{2AMA}
Let $(v_n)$ be a sequence satisfying $v_n\tinf$. 
The object of our first two theorems is to specify the 
logarithmic asymptotic level of the sequences of confidence 
intervals and of confidence bands defined in 
(\ref{Ihat}) and (\ref{Bhat}) respectively.

\begin{theo}
\label{th0}
Assume $(A1)$ holds, set $x\in\bR^d$ such that $f(x)\neq 0$, and assume
that $f$ is twice differentiable at $x$. 
Moreover, assume that $(h_n)$, $(h^*_n)$ and $(v_n)$ satisfy the conditions
\begin{equation}
\label{a}
\frac{nh_n^{*d}}{v_n^2}\rightarrow \infty,\ \ 
v_nh_n^{*2}\rightarrow 0,\ \ \mbox{and}\ \ 
\frac{v_n^2h_n^d}{h_n^{*d}}\rightarrow \infty.
\end{equation}
Then, for any $\de>0$, we have
\[
 \lim_{n\tinf} \frac{v_n^2}{nh_n^{*d}} 
\log \bP \, \left( \,  f(x)\not\in
\left[f_n^*(x)-\delta\frac{\sqrt{f_n(x)\kappa}}{v_n}\ ;\
f_n^*(x)+\delta\frac{\sqrt{f_n(x)\kappa}}{v_n}\right]
\, \right)
 \; = \;  -\frac{\delta^2}{2}.
\]
Moreover, if the additional condition 
$$\frac{nh_n^{*d}}{v_n^2\log(1/h_n^*)}\rightarrow \infty$$
holds, then the sequence of intervals
$$\left[f_n^*(x)-\delta\frac{\sqrt{f_n(x)\kappa}}{v_n}\ ;\
f_n^*(x)+\delta\frac{\sqrt{f_n(x)\kappa}}{v_n}\right]$$
is an asymptotic almost sure sequence of confidence intervals of $f(x)$.
\end{theo}
\medskip

\begin{theo}
\label{th1}
Let $(A1)-(A3)$ hold and assume that $f$ is bounded away from zero on
a compact set $C$. Moreover, assume that $(h_n)$, $(h^*_n)$ and $(v_n)$
satisfy the conditions
\begin{equation}
\label{b}
\frac{nh_n^{*d}}{v_n^2\log(1/h_n^*)}\rightarrow \infty,\ \ 
v_nh_n^{*2}\rightarrow 0,\ \ \mbox{and}\ \ 
\frac{v_n^2h_n^d}{h_n^{*d}}\rightarrow \infty.
\end{equation}
Then, for any $\de>0$, we have
\[
 \lim_{n\tinf} \frac{v_n^2}{nh_n^{*d}} 
\log \bP \, \left( \,  \exists x\in C,\ \ f(x)\not\in
\left[f_n^*(x)-\delta\frac{\sqrt{f_n(x)\kappa}}{v_n}\ ;\
f_n^*(x)+\delta\frac{\sqrt{f_n(x)\kappa}}{v_n}\right]
\, \right)
 \; = \;  -\frac{\delta^2}{2}.
\]
Moreover, the sequence of sets of functions  
$$D_n=\left\{g:\bR^d\rightarrow\bR,\ |g(x)-f_n^*(x)|\leq
\delta\frac{\sqrt{f_n(x)\kappa}}{v_n}\ \forall x\in C\right\}$$
is an asymptotic almost sure sequence of confidence regions of $f$.
\end{theo}
\medskip

\paragraph{Comments on Theorems \ref{th0} and \ref{th1}} 
\begin{itemize}
\item [1)]
The regularity assumption on $f$ in Theorem \ref{th0}
is usually required to establish a CLT for $f_n^*(x)$ in 
the case when $f_n^*$ is defined with a two-order kernel; 
the assumptions on $f$ in Theorem 
\ref{th1} are weaker than those required by 
Bickel and Rosenblatt (1973) to establish the asymptotic 
distribution of the normalized uniform error of the kernel density 
estimator.
\item [2)]
The condition $v_nh_n^{*2}\rightarrow 0$ 
(together with the regularity assumption on $f$) ensures 
that the bias of $f_n^*$ does not interfer. 
In the case $f$ is only once differentiable on $\bR^d$, 
this condition must be replaced by $v_nh_n^*\rightarrow 0$ for Theorems \ref{th0} and \ref{th1} hold.
\item [3)]
The main tool used to prove Theorem \ref{th0} is pointwise MDP 
established for the normalized error of the kernel density estimator, 
whereas the demonstration of Theorem \ref{th1} relies on the use of 
uniform MDP;
that is the reason why the conditions (\ref{b})
in Theorem \ref{th1} are slightly stronger than the conditions (\ref{a})
of Theorem \ref{th0}. Now, as soon as
the conditions of Theorem \ref{th1} hold, pointwise and uniform MDP give 
exactly the same asymptotics. 
This unity in pointwise and uniform MDP differs
from the gap there exists between the nature of the asymptotic 
distribution of the normalized pointwise error of the kernel 
density estimator (given by the CLT) on the one hand, 
and of the normalized uniform error of the kernel 
density estimator (given by Bickel and Rosenblatt (1973))
on the other hand.
\end{itemize}

\subsection{Confidence regions with truncating}
\label{3AMA}
The next theorem allows to set up confidence bands for $f$ on
compact sets $C$ on which $f$ may take the value zero. 
Let the positive real-valued sequences $(h_n)$, $(h^*_n)$,
$(v_n)$ and $(\e_n)$
satisfy the following conditions
\begin{equation}
\label{assumption}
\left\{\begin{array}{l}
{\m{\e_n\ra 0,\ \ \frac{h_n^{*}}{{\e_n}}\rightarrow 0,\ \ 
\frac{h_n^{2}}{\e_n}\rightarrow 0,
}}\\
\\
{\m{v_n\e_n^{3/2}\rightarrow \infty,\ \ 
\frac{v_nh_n^{*2}}{\sqrt{\e_n}}\rightarrow 0,\ \ 
\frac{nh_n^{*d}}{v_n^2\log(1/h_n^*)}\rightarrow \infty,\ \ 
\frac{v_n^2h_n^d\e_n^2}{h_n^{*d}}\rightarrow \infty,}}
\end{array}
\right.
\end{equation}
and let $\tilde T_n$ be the function defined by
(\ref{tildtdef}).

\begin{theo}
\label{th2}
Let $(A1)-(A3)$ hold, assume that there exists $x\in C$ such that
$f(x)\neq 0$, and that the sequences $(h_n)$, $(h^*_n)$,
$(v_n)$ and $(\e_n)$ satisfy (\ref{assumption}).
Then, the conclusions of Theorem \ref{th1} still hold when 
$f_n(x)$ is replaced by $\tilde T_n(x)$.
\end{theo}
\medskip

Let $T_n$ be the function defined by (\ref{Tdef}) with 
$\e_n\leq 1$; 
with the help of Theorem \ref{th2}, we will prove the following result.

\begin{cor}
\label{co1}
Let $(A1)-(A3)$ hold, assume that there exists $x\in C$ such that
$f(x)\neq 0$, and that $(h_n)$, $(h^*_n)$,
$(v_n)$, and $(\e_n)$ satisfy (\ref{assumption}). 
Then, the conclusions of Theorem \ref{th1} still hold when 
$f_n(x)$ is replaced by $T_n(x)$.
\end{cor}
\medskip

The extension of Theorem \ref{th2} and Corollary \ref{co1} to the case
$C=\bR^d$ holds under the additional assumption $(A4)$. 

\begin{theo}
\label{th3}
Let $(A1)-(A4)$ hold, and assume $(h_n)$, $(h^*_n)$,
$(v_n)$ and $(\e_n)$ satisfy (\ref{assumption}). 
Then, the conclusions of Theorem \ref{th1} still hold when 
$f_n(x)$ is replaced by $\tilde T_n(x)$ and $C$ by $\bR^d$.
\end{theo}
\medskip

Let $T_n$ be defined by (\ref{Tdef}) with $\e_n\leq 1$
and $C=\bR^d$.

\begin{cor}
\label{co2}
Let $(A1)-(A4)$ hold, and assume $(h_n)$, $(h^*_n)$,
$(v_n)$, and $(\e_n)$ satisfy (\ref{assumption}). 
Then, the conclusions of Theorem \ref{th1} still hold when 
$f_n(x)$ is replaced by $T_n(x)$ and $C$ by $\bR^d$.
\end{cor}

\paragraph{Remark} Let us mention that Corollaries 
\ref{co1} and \ref{co2} also hold when the sequence $(\e_n)$ is 
constant ($\e_n =\e\in]0,1]$ for all $n$); in the case $\e_n=1$, the width of the 
confidence bands does not depend on the point $x\in C$.

\section{Particular cases} 

In this section, we first give a pointwise analysis of the confidence 
bands provided by Bickel and Rosenblatt (1973), and show how the
MDP approach allows to modify these confidence bands in order to 
obtain confidence bands whose width is of 
the same order as the one of Bickel and Rosenblatt's confidence bands, 
but whose asymptotic level equals one instead of $1-\al<1$; in 
particular, we explicit the choices of the parameters $(h_n^*)$ and 
$(\e_n)$, which give the best convergence rate to one 
of the level of these modified confidence bands. 
Then, we consider the problem of constructing 
confidence bands, which are thinner than those provided by 
Bickel and Rosenblatt, but whose level converges to one slower than 
the level of the modified Bickel and Rosenblatt's 
confidence bands does. 
We give two possible choices of $(h_n^*)$, which both 
correspond to the case our confidence bands are centered at an optimal 
kernel estimator of $f$; for the first choice, the optimality is 
according to the $L^2$ criterion, and, for the second one, to the 
$L^\infty$ criterion.

\subsection{On Bickel and Rosenblatt's confidence bands}

\subsubsection{Pointwise analysis of Bickel and Rosenblatt's 
confidence bands} \label{secBR}

Set $d=1$ and let $C=[c_1,c_2]$ be a bounded interval of $\bR$. 
Bickel and Rosenblatt (1973) construct confidence bands for 
$f$ on $C$ with asymptotic level $1-\al\in]0,1[$ in the case when:
$(i)$ $f$ is bounded away from zero on $C$;
$(ii)$ the kernel $K$ is chosen absolutely continuous on $\bR$ and
such that $\int_\bR K'^2(t)dt\neq 0$; 
$(iii)$ the bandwidth $h_n$ used for the computation of $f_n$
is chosen equal to $(n^{-a})$ with $a\in]\frac{1}{5},\frac{1}{2}[$. 
Their confidence bands are constructed as follows.\\

Set $\al\in]0,1[$, $z_\al$ such that
$\exp(-2\exp(-z_\al))=1-\al$,
and, for all $x\in C$, 
\begin{eqnarray}
I_n(x) & = &
\left[f_n(x)-\frac{\sqrt{f_n(x)\kappa}}{\sqrt{nh_n}}
\sqrt{\log(1/h_n)}\left(\sqrt{2}+u_n+\frac{z_\al}{\sqrt{2}\log(1/h_n)}\right)
\ ;\ \right. \nonumber\\
& & \mbox{} \left.
f_n(x)+\frac{\sqrt{f_n(x)\kappa}}{\sqrt{nh_n}}
\sqrt{\log(1/h_n)}\left(\sqrt{2}+u_n+\frac{z_\al}{\sqrt{2}\log(1/h_n)}
\right)\right]
\label{In def}
\end{eqnarray}
with  
\begin{equation}
\label{iu}
u_n=\frac{1}{\sqrt{2}\log(1/h_n)}\left\{\log\left[\frac{1}{2\pi}
\sqrt{\frac{\int_\bR K'^2(t)dt}{\kappa}}\right]+\log\left[
\frac{c_2-c_1}{\pi}\right]\right\}. 
\end{equation}
Bickel and Rosenblatt (1973) prove that
$$
B_{n,\al} =\cup_{x\in C}I_n(x)
$$
is then a confidence band for $f$ on $C$ with asymptotic
level $1-\al$. 
A straightforward application of the CLT ensures that, 
for each $x\in C$, the asymptotic level of $I_n(x)$ equals one.
Now, Theorem \ref{th0} allows to specify the convergence rate of 
the asymptotic level of the confidence intervals $I_n(x)$ 
toward one. More precisely, the application of Theorem \ref{th0}
with $(h_n)\equiv(h^*_n)$ and $(v_n)\equiv(\sqrt{nh_n/\log(1/h_n)})$,
together with a continuity argument, ensure
that {\bf the logarithmic asymptotic level of each confidence 
interval $I_n(x)$
is $1$ with speed ${\log(1/h_n)}$}. \\

The difference between the asymptotic level of the confidence band 
$B_{n,\al}$ and the asymptotic levels of all the confidence 
intervals $I_n(x)$ is explained by the difference 
between the asymptotic weak behaviour of the uniform error 
of $f_n$ (given by Bickel and 
Rosenblatt's result)  
and the asymptotic weak behaviour of the pointwise error 
of $f_n$ (given by the central limit theorem). Adopting the MDP point
of view, it is note-worthy that this phenomenon corresponds to the
case pointwise MDP hold, but not uniform MDP. As a matter of fact, 
when $d=1$, the sequence $(v_n)\equiv(\sqrt{nh_n/\log(1/h_n)})$ 
fulfills the conditions 
(\ref{a}) required by Theorem \ref{th0}, but not the slightly stronger
conditions (\ref{b}) imposed by  Theorem \ref{th1}.\\

\subsubsection{Improvement of Bickel and Rosenblatt's confidence bands}

The aim of this section is to show how the MDP approach 
allows to improve the confidence bands provided by Bickel and 
Rosenblatt (1973). In a first part, we introduce a translation, which 
allows to provide confidence bands that have the same width as 
the confidence bands provided by Bickel and Rosenblatt, but which 
have a better asymptotic level. In a second part, 
we give a simplification of these translated confidence bands, which 
does affect neither their width order, nor their logarithmic 
asymptotic level. 
Then, we show how we can get rid of the 
condition $f$ is bounded away from zero on $C$. Finally, we 
give the extension to the multivariate framework.

\paragraph{Confidence bands translation}
We consider here Bickel and 
Rosenblatt's framework, that is, the case $d=1$, $C=[c_1,c_2]$, and 
$f$ is bounded away from zero on $C$.\\

Set $(h_n)\equiv (n^{-a})$ with $a\in]\frac{1}{5},\frac{1}{2}[$,  
let $f_n$, $u_n$ and $z_\al$ be defined in the same way as 
in Section \ref{secBR}, and 
$(h_n^*)$ be a bandwidth satisfying the conditions
\begin{equation}
\label{a voir}
{n^a h_n^*}\tinf\ \ \mbox{and}\ \ 
\frac{n^{1-a}h_n^{*4}}{\log n}\ra 0.
\end{equation}
Moreover, let $f^*_n$ be the kernel estimator of $f$ defined with the 
bandwidth $h^*_n$, and, for each $x\in C$, set 
\begin{eqnarray*}
{I}^*_n(x) & = &
\left[{f}^*_n(x)-\frac{\sqrt{f_n(x)\kappa}}{\sqrt{nh_n}}
\sqrt{\log(1/h_n)}\left(\sqrt{2}+u_n+\frac{z_\al}{\sqrt{2}\log(1/h_n)}\right)
\ ;\ \right.\\
& & \mbox{} \left.
{f}^*_n(x)+\frac{\sqrt{f_n(x)\kappa}}{\sqrt{nh_n}}
\sqrt{\log(1/h_n)}\left(\sqrt{2}+u_n+\frac{z_\al}{\sqrt{2}\log(1/h_n)}
\right)\right].
\end{eqnarray*}
Note that, for each $x$ in $C$, ${I}^*_n(x)$ is  
the translation of the confidence interval $I_n(x)$ (defined in 
(\ref{In def})) from the quantity $f_n^*(x)-f_n(x)$. \\

The application of Theorem \ref{th0} (with $d=1$ and 
$(v_n)\equiv(\sqrt{nh_n/\log(1/h_n)})$), 
together with a continuity argument, ensure that, 
{\bf for each $x$ in $C$, 
the logarithmic asymptotic level of ${I}^*_{n}(x)$ is 
equal to $1$ with speed ${h}^*_n\log(1/h_n)/h_n$}. Let us underline 
that the speed obtained for $I^*_{n}(x)$ is faster than the speed 
obtained for $I_{n}(x)$; in other words, the levels of the translated 
intervals ${I}^*_{n}(x)$ go to one faster than the levels of the 
intervals $I_n(x)$. This is explained by the fact that the translated 
intervals ${I}^*_{n}(x)$ are centered at the point ${f}^*_{n}(x)$ rather 
than at the point ${f}_{n}(x)$, and, in view of the conditions 
(\ref{a voir}), the estimator ${f}^*_{n}(x)$ converges to $f(x)$ 
faster than the estimator $f_{n}(x)$ does.
\\

Now, set 
$${B}^*_{n} =\cup_{x\in C}{I}_n^*(x).$$
The application of Theorem \ref{th1} (together with a continuity argument) 
ensures that {\bf ${B}^*_{n}$ is a confidence band for $f$ on $C$
whose logarithmic asymptotic level equals $1$ with speed 
${h}^*_n\log(1/h_n)/h_n$.} \\

The confidence band ${B}^*_{n}$, which is just the translation 
of $B_{n,\al}$ from the quantity $f_n^*-f_n$, has thus the following 
advantages:
\begin{itemize}
\item It has the same width, at each point $x\in C$, as the confidence 
band $B_{n,\al}$ provided by Bickel and Rosenblatt.
\item Its asymptotic level is one instead of being $1-\al<1$.
\item The logarithmic asymptotic level of ${B}^*_{n}$ is the same 
as the logarithmic asymptotic levels of all the intervals ${I}^*_{n}(x)$ 
whose union gives ${B}^*_{n}$, and the intervals ${I}^*_{n}(x)$ 
themselves have a better logarithmic asymptotic level than the 
intervals ${I}_{n}(x)$ whose union gives ${B}_{n,\al}$.
\\
\end{itemize}

Let us mention that the advisable choice of the bandwidth $(h_n^*)$ is 
$(h_n^*)\equiv(n^{-(1-a)/4})$, where $a$ is the parameter which 
defines $(h_n)$. As a matter of fact, among the 
sequences $(h_n^*)\equiv(n^{-b})$ which satisfy (\ref{a voir}), it is the choice 
that maximizes the speed ${h}^*_n\log(1/h_n)/h_n$ (which then 
equals $n^{(5a-1)/4}\log n$). Let us underline that the condition $a>1/5$ 
implies $(1-a)/4<1/5$. 
For this optimal choice of the bandwidth $(h_n^*)$, 
the confidence band ${B}^*_{n}$ (respectively the confidence interval ${I}^*_{n}(x)$) 
is thus centered at an estimator $f_n^*$ (respectively $f_n^*(x)$) whose 
convergence rate is given by the convergence rate of its bias term. Consequently, 
${B}^*_{n}$ (respectively ${I}^*_{n}(x)$) cannot be compared with a
confidence band (respectively confidence interval) centered at $f_n^*$ (respectively $f_n^*(x)$) and provided 
by Bickel and Rosenblatt's result (respectively by the central limit theorem).
The surprising aspect of this result 
is that this optimal choice of $(h_n^*)$ depends on the choice of 
the bandwidth $(h_n)$ and is never the choice for which the estimator 
$f_n^*$ converges at the optimal rate.

\paragraph{Simplification of the translated confidence bands}
The parameters $u_n$ (which depends on the length of the interval 
$C$) and $z_\al$ (which depends on the asymptotic level $\al$), 
which appear in the definitions of the 
intervals $I_n(x)$ and $I_n^*(x)$, play a crucial role in 
Bickel and Rosenblatt's approach. However, they do not have any effect 
in the MDP approach. That is the reason why we 
propose here a simplification of the definition of the confidence 
band $B_n^*$. More precisely, we set 
$$B^{**}_n=\cup_{x\in C}I_n^{**}(x),$$
where, for each $x\in C$,
\begin{eqnarray*}
I_n^{**}(x) & = &
\left[{f}^*_n(x)-\frac{\sqrt{f_n(x)\kappa}}{\sqrt{nh_n}}
\sqrt{\log(1/h_n)}\sqrt{2}
\ ;\ 
{f}^*_n(x)+\frac{\sqrt{f_n(x)\kappa}}{\sqrt{nh_n}}
\sqrt{\log(1/h_n)}\sqrt{2}\right].
\end{eqnarray*}
A straightforward application of Theorems \ref{th0} and \ref{th1} 
ensures that {\bf the logarithmic asymptotic levels of 
${I}^{**}_{n}(x)$ for all $x\in C$ and of ${B}^{**}_{n}$ 
equal $1$ with speed ${h}^*_n\log(1/h_n)/h_n$.}  
In particular, we see that this simplification does not affect 
the logarithmic asymptotic level.

\paragraph{Confidence bands truncating}

Although the simplified confidence band $B^{**}_n$ 
seems very convenient to use, it suffers from the same drawback
as the confidence band $B_{n,\al}$ proposed by Bickel and 
Rosenblatt (1973): its use is conditionned to the fact that the 
density $f$ is bounded away from zero on $C$. In order 
to allow the construction of confidence bands for 
$f$ on intervals $C$ on which $f$ may take the value zero, we now 
introduce the truncated confidence band defined as:
$$B^{***}_n=\cup_{x\in C}I_n^{***}(x),$$
where, for each $x\in C$,
\begin{eqnarray*}
I_n^{***}(x) & = &
\left[{f}^*_n(x)-\frac{\sqrt{T_n(x)\kappa}}{\sqrt{nh_n}}
\sqrt{\log(1/h_n)}\sqrt{2}
\ ;\ 
{f}^*_n(x)+\frac{\sqrt{T_n(x)\kappa}}{\sqrt{nh_n}}
\sqrt{\log(1/h_n)}\sqrt{2}\right],
\end{eqnarray*}
$T_n$ being the truncating function defined in (\ref{Tdef}). 
A straightforward application of Corollary \ref{co1} (respectively of 
Corollary \ref{co2}) in the case $C$ is a compact set (respectively 
in the case $C=\bR$) ensures that, 
if $(\e_n)\equiv(\log n)^{-e}$ with $e\in]0,1[$, then 
{\bf the logarithmic asymptotic level 
of $B^{***}_n$ is $1$ with speed ${h}^*_n\log(1/h_n)/h_n$.} 
In other words, the truncating operation, which allows the construction 
of confidence bands for the density on compact sets on 
which $f$ vanishes or on the whole line, does not affect the 
logarithmic asymptotic level.\\

Let us underline that the advantage of truncating is not only to 
enable the construction of confidence bands for $f$ on 
intervals on which $f$ may take the value zero. Even in 
the case $f$ is bounded away from zero on $C$,
truncating gives, in practice, much better results as soon 
as the length of the interval $C$ is large.\\

\paragraph{The multivariate framework}
As mentionned in the introduction, 
the problem of constructing confidence bands when the probability 
density $f$ is defined on $\bR^d$ has been 
considered by Rosenblatt
(1976). His approach consists in an extension to the $d$-dimensional case
($d>1$) of the results obtained by Bickel and Rosenblatt (1973).
However, in order to enable the construction of confidence bands for $f$
on a compact set $C$ on $\bR^d$ (on which $f$ is bounded away from zero),
Rosenblatt (1976) requires the use of kernels of order $k>d(d+2)/2$, and,
consequently, imposes rather stringent conditions on $f$. 
On the opposite, with the MDP approach, all the 
dimensions are dealt with in the same way. More precisely, let 
the bandwidths $(h_n)$ and $(h_n^*)$ be defined as 
$(h_n)\equiv(n^{-a})$ with $a\in]\frac{1}{d+4},\frac{d+4}{d(d+8)}[$ and 
$(h_n^*)\equiv(n^{-(1-ad)/4})$, the sequence 
$(\e_n)$ as $(\e_n)\equiv(\log n)^{-e}$ with $e\in]0,1[$, 
and set, for each $x\in C$,
\begin{eqnarray*}
I_n^{***}(x) & = &
\left[{f}^*_n(x)-\frac{\sqrt{T_n(x)\kappa}}{\sqrt{nh_n^d}}
\sqrt{\log(1/h_n)}\sqrt{2}
\ ;\ 
{f}^*_n(x)+\frac{\sqrt{T_n(x)\kappa}}{\sqrt{nh_n^d}}
\sqrt{\log(1/h_n)}\sqrt{2}\right],
\end{eqnarray*} 
where the truncating function $T_n$ is defined in (\ref{Tdef}). 
A straightforward application of Corollary \ref{co1} (respectively of 
Corollary \ref{co2}) in the case $C$ is a compact set (respectively 
in the case $C=\bR^d$) ensures that, 
{\bf the logarithmic 
asymptotic level of the confidence band 
$B^{***}_n=\cup_{x\in C}I_n^{***}(x)$ is $1$ with speed 
${h}^{*d}_n\log(1/h_n)/h_n^d$
.} 
Let us mention that this implies the existence of two positive 
functions $\lambda_1^+$ and $\lambda_1^-$ which go to infinity with a 
logarithmic rate, and such that 
$$
\exp\left(-\frac{\de^2}{2}n^{([d+4]a-1)d/4}
\lambda_1^-(n)\right)\leq 
\bP\left(\exists x\in C,\ f(x)\not\in I_n^{***}(x)\right)\leq
\exp\left(-\frac{\de^2}{2}n^{([d+4]a-1)d/4}
\lambda_1^+(n)\right).
$$

\subsection{Thinner confidence bands} 

The width order of the confidence band $B^{***}_n$ is 
$(\log n)^{1/2}n^{-b}$ 
with $b<2/(d+4)$; this width might seem too large, and thinner 
confidence bands might be prefered, although the convergence rate 
to 1 of their asymptotic level is slower. 

The smallest possible width of confidence bands whose width does not 
depend on $\sqrt{f(x)}$ and whose asymptotic level equals $1-\al<1$ is, according 
to the minimax theory, $M[(\log n)/n]^{2/(d+4)}$ where the constant 
$M$ depends on some known bound of $\|f\|_\infty$ and $\|f''\|_\infty$ 
(see Ibragimov and Hasminskii (1981),
Donoho and Liu (1991), Donoho (1994), and Tsybakov (2004)). 
This optimal width can not be reached in the case the 
width of the confidence bands depends on $\sqrt{f(x)}$; as a matter of 
fact, for a large class of densities (which includes the standard 
Gaussian density), the sequence 
$[n/\log n]^{2/(d+4)}\|(f_n-f)/\sqrt{f}\|_\infty$ is known to be not 
stochastically bounded (see Gin\'e and Guillou (2002), pp. 918).

In this section, we give two examples of choices of the parameters 
$(h_n^*)$, $(v_n)$, $(h_n)$, and $(\e_n)$, which lead to confidence 
bands whose width order is close to $[(\log n)/n]^{2/(d+4)}$.

\begin{itemize}
\item
Set $(h_n^*)\equiv (c^* n^{-1/(d+4)})$ with $c^*>0$; this choice 
corresponds to the case the confidence bands provided by the MDP 
approach are centered at the kernel estimator, which minimizes the  
(integrated) mean squared error. 
For this choice of bandwidth, the sequence $(v_n)$ can be chosen 
as:
$$(v_n)\equiv\left(v^*\frac{n^{2/(d+4)}}{(\log n)^a}\right)
\ \ \mbox{with}\ \ v^*>0 \ \ \mbox{and}\ \ a>\frac{1}{2},$$
the sequence $(h_n)$ can be chosen equal to $(h_n^*)$ 
or to $(c[n/\log n]^{-1/(d+4)})$ with $c>0$, and 
the sequence $(\e_n)$ as:
$$(\e_n)\equiv\left(\e^*(\log n)^{-e}\right)
\ \ \mbox{with}\ \ \e^*>0 \ \ \mbox{and}\ \ e<2a.$$
The application of Corollary \ref{co1} (in the case $C$ is a compact set) 
or of Corollary \ref{co2} (in the case $C=\bR^d$)
ensures that the logarithmic asymptotic level of the 
confidence bands defined as $B_n=\cup_{x\in C}{I}_n(x)$
with 
\begin{equation}
\label{I}
{I}_n(x)=\left[f_n^*(x)-\de
\frac{\sqrt{{T}_n(x)\kappa}}{v_n}\ ;\
f_n^*(x)+\de
\frac{\sqrt{{T}_n(x)\kappa}}{{v_n}}\right]
\end{equation}
is then equal to $\de^2/2$ with speed $nh_n^{*d}/v^2_n$.
Consequently, there exist two positive functions $\lambda_2^+$ 
and $\lambda_2^-$ which go to infinity with a 
logarithmic rate, and such that 
$$
n^{-\frac{\de^2}{2}\lambda_2^-(n)}
\leq 
\bP\left(\exists x\in C,\ f(x)\not\in I_n(x)\right)\leq
n^{-\frac{\de^2}{2}\lambda_2^+(n)}.
$$

\item
Set $(h_n^*)\equiv (c^* [n/\log n]^{-1/(d+4)})$ with $c^*>0$; 
this choice 
corresponds to the case the confidence bands 
are centered at the kernel estimator, which minimizes the  
uniform error. 
For this choice of bandwidth, we can construct confidence bands whose width is arbitrarily close to $[(\log n)/n]^{2/(d+4)}$, by choosing the sequence $(v_n)$ as
$$(v_n)\equiv\left(v^*\left[\frac{n}{\log n}\right]^{2/(d+4)}
\frac{1}{(\log\log n)^a}\right)
\ \ \mbox{with}\ \ v^*>0 \ \ \mbox{and}\ \ a>0,$$
the sequence $(h_n)$ equal to $(h_n^*)$, and 
the sequence $(\e_n)$ as
$$(\e_n)\equiv\left(\e^*(\log\log n)^{-e}\right)
\ \ \mbox{with}\ \ \e^*>0 \ \ \mbox{and}\ \ e<2a.$$
The application of Corollaries \ref{co1} and \ref{co2} 
ensures that logarithmic asymptotic level of the 
confidence bands defined as $B_n=\cup_{x\in C}{I}_n(x)$
with ${I}_n(x)$ defined in (\ref{I}) 
is then equal to $\de^2/2$ 
with speed $nh_n^{*d}/v^2_n$.
Accordingly, 
$$
(\log n)^{-\frac{\de^2}{2}\lambda_3^-(n)}
\leq 
\bP\left(\exists x\in C,\ f(x)\not\in I_n(x)\right)\leq
(\log n)^{-\frac{\de^2}{2}\lambda_3^+(n)}
$$
where $\lambda_3^+$ and $\lambda_3^-$ are two positive functions, 
which go to infinity with a rate in $\log\log$.
\end{itemize}

Let us finally mention that, in all the previous examples, the truncating 
function $T_n$ can be replaced by the function $\ti T_n$ defined in 
(\ref{tildtdef}).

\section{Proofs} 

We first give a unified proof for all the almost sure parts of our results in Section \ref{Sec Unif}. Then, 
Theorem \ref{th0} is proved in Section \ref{Sec thm 1},
Theorems \ref{th1}, \ref{th2}, and \ref{th3} 
in Section \ref{Sec thm 123}, and 
Corollaries \ref{co1} and \ref{co2} in Section \ref{3.4}.

\subsection{Unified proof for all the almost sure parts of our results} 
\label{Sec Unif}

The proof relies on the use of both conditions 
$$v_nh_n^{*2}\rightarrow 0 \ \ \mbox{and}\ \
\frac{nh_n^{*d}}{v_n^2\log(1/h_n^*)}\rightarrow \infty.$$
Set $\g= \de^2/2$, $w_n={nh_n^{*d}}/v_n^2$, 
$\rho\in]0,1/(d+4)[$, and $M>1/\rho$. On the one hand, the condition
$v_nh_n^{*2}\rightarrow 0$ implies that, for $n$ large enough,
$$\exp(-\g w_n/2)\leq \exp\left[\frac{-\g nh_n^{*(d+4)}}{2}\right].$$
On the other hand, the condition 
$nh_n^{*d}/[v_n^2\log(1/h_n^*)]\tinf$ implies that, for $n$ large enough,
$nh_n^{*d}/v_n^2\geq 2M\log(1/h_n^*)/\g$, and thus
$$\exp(-\g w_n/2)\leq h_n^{*M}.$$
It follows that 
$$\exp(-\g w_n/2)\leq
\left\{
\begin{array}{ll}
{\m \exp\left(-\g n^{1-(d+4)\rho}/2\right)} &
\mbox{if $h_n^*\geq n^{-\rho}$} \\
{\m n^{-M\rho}} &
\mbox{if $h_n^*\leq n^{-\rho}$},
\end{array} \right. $$
and thus $\sum_n \exp(-\g w_n/2)<\infty$. 
The almost sure parts of our results then follow from the application of 
Proposition \ref{aspro}.

\subsection{Proof of Theorem \ref{th0}}
\label{Sec thm 1}

Set $\de>0$ and $\eta\in]0,1[$. On the one hand, we have
\begin{eqnarray*}
\lefteqn{\bP \, \left( \,  f(x)\not\in
\left[f_n^*(x)-\delta\frac{\sqrt{f_n(x)\kappa}}{v_n}\ ;\
f_n^*(x)+\delta\frac{\sqrt{f_n(x)\kappa}}{v_n}\right]
\, \right)}\\
& \leq &
\bP \left[ {v_n|f_n^*(x)-f(x)|}
 > \delta\sqrt{f_n(x)\kappa}   \ \ \mbox{and}\ \ \frac{f_n(x)}{f(x)}>
\frac{1}{1+\eta}\right]
+\bP \left[ \frac{f_n(x)}{f(x)}\leq 
\frac{1}{1+\eta} \right] \\ 
& \leq &
\bP \left[ {v_n|f_n^*(x)-f(x)|}
 > \frac{\delta{\sqrt{f(x)\kappa}}}{\sqrt{1+\eta}}  \right]
+\bP \left[ f(x)-f_n(x) \geq \frac{\eta f(x)}{1+\eta}  \right].
\end{eqnarray*}
Now, Theorem 4 in Mokkadem, Pelletier 
and Worms (2005) ensures that
$$\lim_{n\tinf} \frac{v_n^2}{nh_n^{*d}} 
\log \bP \left[ {v_n|f_n^*(x)-f(x)|}
 > \frac{\delta{\sqrt{f(x)\kappa}}}{\sqrt{1+\eta}}  \right]
\; = \; \frac{-\delta^2}{2(1+\eta)},$$
and, since ${v_n^2h_n^d}/{h_n^{*d}}\rightarrow\infty$, the 
application of Corollary 1 in Mokkadem, Pelletier 
and Worms (2005) gives 
\begin{eqnarray*} 
\limsup_{n\tinf} \frac{v_n^2}{nh_n^{*d}} 
\log\bP \left[ f(x)-f_n(x) \geq \frac{\eta f(x)}{1+\eta} \right]   
& \leq &
\limsup_{n\tinf} \frac{v_n^2h_n^d}{h_n^{*d}} \left[\frac{1}{nh_n^{d}} 
\log\bP \left[ |f(x)-f_n(x)| \geq \frac{\eta f(x)}{1+\eta}  \right] \right].\\
& = & -\infty.
\end{eqnarray*}
We thus deduce that
\begin{equation}
\label{1etth1}
\limsup_{n\tinf} \frac{v_n^2}{nh_n^{*d}} 
\log 
{\bP \, \left( \,  f(x)\not\in
\left[f_n^*(x)-\delta\frac{\sqrt{f_n(x)\kappa}}{v_n}\ ;\
f_n^*(x)+\delta\frac{\sqrt{f_n(x)\kappa}}{v_n}\right]
\, \right)}  
\; \leq \; \frac{-\delta^2}{2(1+\eta)}.
\end{equation}
On the other hand, we note that
\begin{eqnarray*}
\lefteqn{\bP \, \left( \,  f(x)\not\in
\left[f_n^*(x)-\delta\frac{\sqrt{f_n(x)\kappa}}{v_n}\ ;\
f_n^*(x)+\delta\frac{\sqrt{f_n(x)\kappa}}{v_n}\right]
\, \right)}\\
& \geq &
\bP \left[ {v_n|f_n^*(x)-f(x)|}
 > \delta {\sqrt{f(x)\kappa}}\sqrt{\frac{f_n(x)}{f(x)}}    \ \ \mbox{and}\ \ 
\sqrt{\frac{f_n(x)}{f(x)}}\leq\sqrt{1+\eta} \right]  \\
& \geq &
\bP \left[ {v_n|f_n^*(x)-f(x)|}
 > \delta {\sqrt{(1+\eta)f(x)\kappa}}    \ \ \mbox{and}\ \ 
\sqrt{\frac{f_n(x)}{f(x)}}\leq\sqrt{1+\eta} \right]  \\
& \geq &
\bP \left[ {v_n|f_n^*(x)-f(x)|}
 > {\delta} {\sqrt{(1+\eta)f(x)\kappa}}   \right]  
-\bP \left[ f_n(x) > (1+\eta) f(x)  \right]
\\ 
& \geq &
\bP \left[ {v_n|f_n^*(x)-f(x)|}
 > {\delta} {\sqrt{(1+\eta)f(x)\kappa}}   \right]
-\bP \left[ |f_n(x)-f(x)| > \eta f(x)  \right],
\end{eqnarray*}
and the application of Corollary 1 and Theorem 4 in Mokkadem, Pelletier 
and Worms (2005) leads to 
\begin{equation}
\label{2etth1}
\liminf_{n\tinf} \frac{v_n^2}{nh_n^{*d}} 
\log 
{\bP \, \left( \,  f(x)\not\in
\left[f_n^*(x)-\delta\frac{\sqrt{f_n(x)\kappa}}{v_n}\ ;\
f_n^*(x)+\delta\frac{\sqrt{f_n(x)\kappa}}{v_n}\right]
\, \right)}   
\; \geq \; \frac{-\delta^2(1+\eta)}{2}.
\end{equation}
Since $\eta$ can be taken arbitrarily close to zero, Theorem \ref{th0}
is a straightforward consequence of (\ref{1etth1}) and (\ref{2etth1}).

\subsection{Proof of Theorems \ref{th1}, \ref{th2}, and 
\ref{th3}}
\label{Sec thm 123}

The proof of Theorems \ref{th1}, \ref{th2}, and \ref{th3} 
will require the application of Lemmas 
\ref{essai 1 lprobadusup} and \ref{essai 2 lprobadusup} 
below. We first state these Lemmas, whose proof is 
postponed in the appendix (see Section \ref{Sec App 1}). 
Then, we prove Theorems \ref{th3}, \ref{th2}, 
and \ref{th1} in Sections \ref{3.3.1}, \ref{3.3.2}, 
and \ref{3.3.3} respectively.\\

Let $C_n$ be a sequence of compact sets of $\bR^d$ and set 
$w_n=\sup\{\|x\|,\ x\in C_n\}$. 
Moreover, set $\xi\in]0,1[$, $\zeta>1$, and 
\begin{eqnarray}
{\ti U}_n & = & \left\{x\in C_n,\ f_n(x)\geq\e_n\right\}, 
\label{2 global}\\
{\ti W}_n & = & \left\{x\in C_n,\ f_n(x)<\e_n\right\}, 
\label{3 global}\\
U_n(\xi) & = & \left\{x\in C_n,\ f(x)\geq\xi\e_n\right\},  
\label{4 global}\\
W_n(\zeta) & = & \left\{x\in C_n,\ f(x)\leq\zeta\e_n\right\}.
\label{5 global}
\end{eqnarray}

\begin{lemma}
\label{essai 1 lprobadusup} 
Assume that (A1)-(A3) hold, and that $(v_n)$, $(h_n)$, 
$(h_n^*)$, $(\e_n)$ satisfy (\ref{assumption}). 
Moreover, assume that $(w_n)$ fulfills the condition
\begin{equation}
\label{wdef}
\lim_{n\tinf}\frac{v_n^2\log w_n}{nh_n^{*d}} =0.
\end{equation}
Then, for any $\delta>0$, 
$$
 \limsup_{n\tinf} \frac{v_n^2}{nh_n^{*d}} 
\log \bP \, \left[ \,  \sup_{x\in U_n(\xi)}
\frac{v_n|f_n^*(x)-f(x)|}{\sqrt{f_n(x)\ka}}
 \geq \delta  
\ \ \mbox{and}\ \ 
\inf_{x\in U_n(\xi)}f_n(x)>0
\, \right] \; \leq \;  -\frac{\delta^2}{2}.
$$
\end{lemma}
\vspace{0.5cm}

\begin{lemma}
\label{essai 2 lprobadusup} 
Assume that (A1)-(A3) hold, and that $(v_n)$, $(h_n)$, 
$(h_n^*)$, $(\e_n)$ and $(w_n)$ satisfy (\ref{assumption}) 
and (\ref{wdef}).
For any $\delta>0$, 
$$
 \limsup_{n\tinf} \frac{v_n^2}{nh_n^{*d}} 
\log \bP \, \left[ \,  \sup_{x\in W_n(\zeta)}
\frac{v_n|f_n^*(x)-f(x)|}{\sqrt{\e_n\kappa}}
 \geq \delta  \, \right] \; \leq \; -\delta^2\left(1-\frac{\zeta}{2}\right).
$$
\end{lemma}
\vspace{0.5cm}

\subsubsection{Proof of Theorem \ref{th3}}
\label{3.3.1} 

To prove Theorem \ref{th3}, we first establish the upper bound 
\begin{equation}
\label{Maj pour thm4}
\limsup_{n\tinf} \frac{v_n^2}{nh_n^{*d}} 
\log \bP  \left(   \exists x\in\bR^d,\ \ f(x)\not\in
\left[f_n^*(x)-\delta\frac{\sqrt{\tilde T_n(x)\kappa}}{v_n}\ ;\
f_n^*(x)+\delta\frac{\sqrt{\tilde T_n(x)\kappa}}{v_n}\right]
 \right)
 \leq   -\frac{\delta^2}{2},
\end{equation}
and then prove the lower bound 
\begin{equation}
\label{Min pour thm4}
 \liminf_{n\tinf} \frac{v_n^2}{nh_n^{*d}} 
\log \bP  \left(   \exists x\in\bR^d,\ \ f(x)\not\in
\left[f_n^*(x)-\delta\frac{\sqrt{\tilde T_n(x)\kappa}}{v_n}\ ;\
f_n^*(x)+\delta\frac{\sqrt{\tilde T_n(x)\kappa}}{v_n}\right]
 \right)
 \geq   -\frac{\delta^2}{2}.
\end{equation}
\\
Throughout the proof, we set 
$C_n=\left\{x \in\bR^d,\ \|x\|\leq\e_n^{-2/q}\right\}$
(we thus have $w_n=\e_n^{-2/q}$).

\paragraph{Proof of the upper bound (\ref{Maj pour thm4})}
We have  
\begin{eqnarray*}
& &\bP \, \left( \,  \exists x\in \bR^d,\ f(x)\not\in
\left[f_n^*(x)-\delta\frac{\sqrt{\ti T_n(x)\kappa}}{v_n}\ ;\
f_n^*(x)+\delta\frac{\sqrt{\ti T_n(x)\kappa}}{v_n}\right]
\, \right)  
\\
&  =  & 
\bP \, \left[ \, \sup_{x\in  C_n} 
  \frac{{v_n}\left|f_n^*(x)-f(x)\right|}{\sqrt{\tilde T_n(x)\ka}}
  \, > \,  \delta \, \right] 
+\bP \, \left[ \, \sup_{x\in C_n^c} 
  \frac{{v_n}\left|f_n^*(x)-f(x)\right|}{\sqrt{\tilde T_n(x)\ka}}
  \, > \,  \delta \, \right] \\
& \leq & 
\bP \, \left[ \, \sup_{x\in {\ti U}_n} 
  \frac{{v_n}\left|f_n^*(x)-f(x)\right|}{\sqrt{ f_n(x)\ka}}  \, > \,  
\delta \, \right] 
+\bP \, \left[ \, \sup_{x\in {\ti W}_n} 
  \frac{{v_n}\left|f_n^*(x)-f(x)\right|}{\sqrt{\e_n\ka}}  \, > \,  
\delta \, \right]\\
& & \mbox{~}
+ \bP \, \left[ \, \sup_{x\in C_n^c} 
  \frac{{v_n}\left|f_n^*(x)-f(x)\right|}{\sqrt{\tilde T_n(x)\ka}}
  \, > \,  \delta \, \right] \\
& \leq & 
\bP \, \left[ \, \sup_{x\in {\ti U}_n\cap U_n(\xi)} 
  \frac{{v_n}\left|f_n^*(x)-f(x)\right|}{\sqrt{ f_n(x)\ka}} 
 \, \geq \,  \delta \, \right]
+\bP \, \left[ \, {\ti U}_n\cap [U_n(\xi)]^c\neq\emptyset \, \right]\\
& & \mbox{~}
+ \bP \, \left[ \, \sup_{x\in W_n(\zeta)} 
  \frac{{v_n}\left|f_n^*(x)-f(x)\right|}{\sqrt{\e_n\ka}}  \, \geq \,  
\delta \, \right] 
+\bP \, \left[ \, {\ti W}_n\cap [W_n(\zeta)]^c\neq\emptyset \, \right]\\
& & \mbox{~}
+ \bP \, \left[ \, \sup_{x\in C_n^c} 
  \frac{{v_n}\left|f_n^*(x)-f(x)\right|}{\sqrt{\tilde T_n(x)\ka}}
  \, > \,  \delta \, \right] 
\\
& \leq &
\bP \, \left[ \, \sup_{x\in U_n(\xi)} 
  \frac{{v_n}\left|f_n^*(x)-f(x)\right|}{\sqrt{ f_n(x)\ka}}  \, \geq \,  
\delta 
\ \ \mbox{and}\ \ \inf_{x\in U_n(\xi)}f_n(x)>0
\, \right]
+\bP \, \left[ \, {\ti U}_n\cap [U_n(\xi)]^c\neq\emptyset \, \right]\\
& & \mbox{~}
+ \bP \, \left[ \, \sup_{x\in W_n(\zeta)} 
  \frac{{v_n}\left|f_n^*(x)-f(x)\right|}{\sqrt{\e_n\ka}}  \, \geq \,  
\delta \, \right] 
+\bP \, \left[ \, {\ti W}_n\cap [W_n(\zeta)]^c\neq\emptyset \, \right]\\
& & \mbox{~}
+ \bP \, \left[ \, \sup_{x\in C_n^c} 
  \frac{{v_n}\left|f_n^*(x)-f(x)\right|}{\sqrt{\tilde T_n(x)\ka}}
  \, > \,  \delta \, \right], 
\end{eqnarray*}
so that
\begin{eqnarray}
& &
\limsup_{n\tinf} \frac{v_n^2}{nh_n^{*d}} 
\log \bP  \left(   \exists x\in\bR^d,\ \ f(x)\not\in
\left[f_n^*(x)-\delta\frac{\sqrt{\tilde T_n(x)\kappa}}{v_n}\ ;\
f_n^*(x)+\delta\frac{\sqrt{\tilde T_n(x)\kappa}}{v_n}\right]
 \right) \nonumber \\
& \leq & 
\max\left\{\limsup_{n\tinf} \frac{v_n^2}{nh_n^{*d}} 
\bP \, \left[  \sup_{x\in U_n(\xi)} 
  \frac{{v_n}\left|f_n^*(x)-f(x)\right|}{\sqrt{ f_n(x)\ka}}  
\geq   \delta 
\ \ \mbox{and}\ \ \inf_{x\in U_n(\xi)}f_n(x)>0 \right]\ ; \ 
\right. \nonumber \\ & & \mbox{~} \left.
\limsup_{n\tinf} \frac{v_n^2}{nh_n^{*d}} 
\bP  \left[  {\ti U}_n\cap [U_n(\xi)]^c\neq\emptyset \right]
\ ; \
\limsup_{n\tinf} \frac{v_n^2}{nh_n^{*d}} 
\bP  \left[  \sup_{x\in W_n(\zeta)} 
  \frac{{v_n}\left|f_n^*(x)-f(x)\right|}{\sqrt{\e_n\ka}}  \geq 
\delta  \right] \ ; \
\right. \nonumber \\ & & \mbox{~} \left.
\limsup_{n\tinf} \frac{v_n^2}{nh_n^{*d}} 
\bP  \left[  {\ti W}_n\cap [W_n(\zeta)]^c\neq\emptyset  \right]
\ ; \
\limsup_{n\tinf} \frac{v_n^2}{nh_n^{*d}} 
\bP  \left[  \sup_{x\in C_n^c} 
  \frac{{v_n}\left|f_n^*(x)-f(x)\right|}{\sqrt{\tilde T_n(x)\ka}}
   >  \delta  \right]
\right\}.
\label{Max}
\end{eqnarray}

\begin{itemize}
\item 
The application of Lemma \ref{essai 1 lprobadusup} 
ensures that 
\begin{eqnarray}
& & \limsup_{n\tinf} \frac{v_n^2}{nh_n^{*d}} 
\log \bP \, \left[ \,  \sup_{x\in U_n(\xi)}
\frac{v_n|f_n^*(x)-f(x)|}{\sqrt{f_n(x)\ka}}
 \geq \delta  
\ \ \mbox{and}\ \ 
\inf_{x\in U_n(\xi)}f_n(x)>0
\, \right]  \leq   -\frac{\delta^2}{2} ~~~~
\label{ancien 1part lprobadusup} 
\end{eqnarray}
and the one of Lemma \ref{essai 2 lprobadusup} gives 
\begin{eqnarray}
& &\limsup_{n\tinf} \frac{v_n^2}{nh_n^{*d}} 
\log \bP \, \left[ \,  \sup_{x\in W_n(\zeta)}
\frac{v_n|f_n^*(x)-f(x)|}{\sqrt{\e_n\kappa}}
 \geq \delta  \, \right]  \leq  
-\delta^2\left(1-\frac{\zeta}{2}\right).
\label{ancien 2part lprobadusup}
\end{eqnarray}
\item
The proof of the following upper bound is quite 
technical, and is postponed in the appendix (see 
Section \ref{Sec App 2}).
\begin{eqnarray}
& &\limsup_{n\tinf} \frac{v_n^2}{nh_n^{*d}} 
\log \bP \, \left[ \,  \sup_{x\in C_n^c}
\frac{v_n|f_n^*(x)-f(x)|}{\sqrt{\tilde T_n(x)\kappa}}
 \geq \delta  \, \right]  =  -\infty.
\label{ancien deduitdeapplgg}
\end{eqnarray}
\item 
Since 
\begin{eqnarray*}
\bP \, \left[ \, {\ti U}_n\cap [U_n(\xi)]^c\neq\emptyset \, \right]
& = &
\bP \, \left[ \, \exists x_0\in C_n,\  f_n\left(x_0\right)\geq\e_n\ \mbox{and}
\ f\left(x_0\right)\leq\xi\e_n \, \right]\\
& \leq &
\bP \, \left[ \, \exists x_0\in C_n,\  
f_n\left(x_0\right)-f\left(x_0\right)\geq(1-\xi)\e_n \, \right]\\
& \leq &
\bP \, \left[ \, \frac{1}{\e_n}\sup_{x\in C_n}  
\left|f_n(x)-f(x)\right|\geq(1-\xi) \, \right],
\end{eqnarray*}
we get, by application of Theorem 5 in Mokkadem, Pelletier 
and Worms (2005), 
\begin{eqnarray}
\lefteqn{\limsup_{n\tinf}\frac{v_n^2}{nh_n^{*d}} 
\log \bP \, \left[ \,  {\ti U}_n\cap [U_n(\xi)]^c\neq\emptyset
\, \right]} \nonumber\\
& \leq &
\limsup_{n\tinf}
\frac{v_n^2\e_n^2h_n^d}{h_n^{*d}}\left\{\frac{1}{nh_n^{d}\e_n^2}
\log \bP \, \left[ \, \frac{1}{\e_n}\sup_{x\in C_n}  
\left|f_n(x)-f(x)\right|\geq(1-\xi) \, \right]\right\} \nonumber\\
& = &
-\infty .
\label{a2refeq}
\end{eqnarray}
\item 
Similarly, since 
\begin{eqnarray*}
\bP \, \left[ \, {\ti W}_n\cap [W_n(\zeta)]^c\neq\emptyset \, \right]
& = &
\bP \, \left[ \, \exists x_0\in C_n,\  f_n\left(x_0\right)<\e_n\ \mbox{and}
\ f\left(x_0\right)>\zeta\e_n \, \right]\\
& \leq &
\bP \, \left[ \, \exists x_0\in C_n,\  
f\left(x_0\right)-f_n\left(x_0\right)>(\zeta -1)\e_n \, \right]\\
& \leq &
\bP \, \left[ \, \frac{1}{\e_n}\sup_{x\in C_n}  
\left|f_n(x)-f(x)\right|>(\zeta -1) \, \right],
\end{eqnarray*}
the application of Theorem 5 in Mokkadem, Pelletier 
and Worms (2005) gives
\begin{eqnarray}
\lefteqn{\limsup_{n\tinf}\frac{v_n^2}{nh_n^{*d}} 
\log \bP \, \left[ \,  {\ti W}_n\cap [W_n(\zeta)]^c\neq\emptyset
\, \right]} \nonumber\\ 
& \leq &
\limsup_{n\tinf}
\frac{v_n^2\e_n^2h_n^d}{h_n^{*d}}\left\{\frac{1}{nh_n^{d}\e_n^2}
\log \bP \, \left[ \, \frac{1}{\e_n}\sup_{x\in C_n}  
\left|f_n(x)-f(x)\right|\geq(\zeta -1) \, \right]\right\} \nonumber\\
& = &
-\infty .
\label{bis a2refeq}
\end{eqnarray}
\end{itemize}
The combination of (\ref{Max})-(\ref{bis a2refeq}) leads to  
$$\limsup_{n\tinf} \frac{v_n^2}{nh_n^{*d}} 
\log \bP  \left(   \exists x\in\bR^d,\ \ f(x)\not\in
\left[f_n^*(x)-\delta\frac{\sqrt{\tilde T_n(x)\kappa}}{v_n}\ ;\
f_n^*(x)+\delta\frac{\sqrt{\tilde T_n(x)\kappa}}{v_n}\right]
 \right)
 \leq -\delta^2\left(1 -\frac{\zeta}{2}\right). $$
Since this last upper bound holds for any $\zeta>1$, 
the upper bound (\ref{Maj pour thm4}) follows.

\paragraph{Proof of the lower bound (\ref{Min pour thm4})}

Set $x_0\in \cap_n C_n$ such that $f(x_0)\neq 0$, and set $\eta\in]0,1[$.
Moreover, let $n$ be large enough so that $f(x_0)\geq\e_n/(1-\eta)$.
We then have:
\begin{eqnarray*}
\lefteqn{
{{v_n}\left|f_n^*(x_0)-f(x_0)\right|}
\geq\de{\sqrt{f_n(x_0)\ka}}\ \ \mbox{and}\ \ f_n(x_0)\geq (1-\eta)f(x_0)
} \\
& \imp &
\frac{{v_n}\left|f_n^*(x_0)-f(x_0)\right|}{\sqrt{f_n(x_0)\ka}}
\geq\de\ \ \mbox{and}\ \ f_n(x_0)\geq \e_n  \\
& \imp &
\frac{{v_n}\left|f_n^*(x_0)-f(x_0)\right|}{\sqrt{f_n(x_0)\ka}}
\geq\de\ \ \mbox{and}\ \ x_0\in{\ti U}_n\\
& \imp &
\sup_{x\in {\ti U}_n} 
  \frac{{v_n}\left|f_n^*(x)-f(x)\right|}{\sqrt{\tilde T_n(x)\ka}}\geq\de\\
& \imp &
\sup_{x\in C_n} 
  \frac{{v_n}\left|f_n^*(x)-f(x)\right|}{\sqrt{\tilde T_n(x)\ka}}\geq\de . 
\end{eqnarray*}
It follows that
\begin{eqnarray*}
\lefteqn{
\bP \, \left[ \, \sup_{x\in C_n} 
  \frac{{v_n}\left|f_n^*(x)-f(x)\right|}{\sqrt{\tilde T_n(x)\ka}}
  \, \geq \,  \delta \, \right]
}\\
& \geq &
\bP \, \left[ \,  
  {{v_n}\left|f_n^*(x_0)-f(x_0)\right|}
  \, \geq \,  \delta{\sqrt{f_n(x_0)\ka}} \, \right]
- \bP \, \left[ \,  f_n(x_0)< (1-\eta)f(x_0)
   \, \right]\\
& \geq &
\bP \, \left[ \,  
  {{v_n}\left|f_n^*(x_0)-f(x_0)\right|}
  \, \geq \,  \delta {\sqrt{f_n(x_0)\ka}} \, \right]
- \bP \, \left[ \,  f(x_0)-f_n(x_0)> \eta f(x_0)
   \, \right].
\end{eqnarray*}
Since the application of Corollary 1 in Mokkadem, 
Pelletier and Worms (2005) ensures that
\begin{eqnarray*}
\limsup_{n\tinf} \frac{v_n^2}{nh_n^{*d}} 
\log \bP \, \left[ \,  f(x_0)-f_n(x_0)> 
\eta f(x_0)   \, \right]
& = &
-\infty,
\end{eqnarray*}
the application of Theorem \ref{th0} leads to
\begin{eqnarray}
\liminf_{n\tinf} \frac{v_n^2}{nh_n^{*d}} \log\bP \, \left[ \, \sup_{x\in C_n} 
  \frac{{v_n}\left|f_n^*(x)-f(x)\right|}{\sqrt{\tilde T_n(x)\ka}}
  \, \geq \,  \delta \, \right]
& \geq &
-\frac{\delta^2}{2}.
\label{ancien lminoration}
\end{eqnarray}
Noting that 
\begin{eqnarray*}
&&\bP \, \left( \,  \exists x\in \bR^d,\ f(x)\not\in
\left[f_n^*(x)-\delta\frac{\sqrt{\ti T_n(x)\kappa}}{v_n}\ ;\
f_n^*(x)+\delta\frac{\sqrt{\ti T_n(x)\kappa}}{v_n}\right]
\, \right)\\
&&\geq 
\bP \, \left[ \, \sup_{x\in C_n} 
  \frac{{v_n}\left|f_n^*(x)-f(x)\right|}{\sqrt{\tilde T_n(x)\ka}}
  \, > \,  \delta \, \right],
\end{eqnarray*}
the lower bound (\ref{Min pour thm4}) follows, which 
concludes the proof of Theorem \ref{th3}.

\subsubsection{Proof of Theorem \ref{th2}}
\label{3.3.2}  

Set $C_n=C$ for all $n$, 
set $\xi\in]0,1[$, $\zeta>1$, and let 
${\ti U}_n$, ${\ti W}_n$, $U_n(\xi)$, and 
$W_n(\zeta)$ be defined according to (\ref{2 global}), 
(\ref{3 global}), (\ref{4 global}), 
and (\ref{5 global}) respectively.

\paragraph{Upper bound}
Following the proof of (\ref{Max}), we have: 
\begin{eqnarray*}
& &
\limsup_{n\tinf} \frac{v_n^2}{nh_n^{*d}} 
\log \bP  \left(   \exists x\in C,\ \ f(x)\not\in
\left[f_n^*(x)-\delta\frac{\sqrt{\tilde T_n(x)\kappa}}{v_n}\ ;\
f_n^*(x)+\delta\frac{\sqrt{\tilde T_n(x)\kappa}}{v_n}\right]
 \right) \nonumber \\
& \leq & 
\max\left\{\limsup_{n\tinf} \frac{v_n^2}{nh_n^{*d}} 
\bP \, \left[  \sup_{x\in U_n(\xi)} 
  \frac{{v_n}\left|f_n^*(x)-f(x)\right|}{\sqrt{ f_n(x)\ka}}  
\geq   \delta 
\ \ \mbox{and}\ \ \inf_{x\in U_n(\xi)}f_n(x)>0 \right]\ ; \ 
\right. \nonumber \\ & & \mbox{~} \left.
\limsup_{n\tinf} \frac{v_n^2}{nh_n^{*d}} 
\bP  \left[  {\ti U}_n\cap [U_n(\xi)]^c\neq\emptyset \right]
\ ; \
\limsup_{n\tinf} \frac{v_n^2}{nh_n^{*d}} 
\bP  \left[  \sup_{x\in W_n(\zeta)} 
  \frac{{v_n}\left|f_n^*(x)-f(x)\right|}{\sqrt{\e_n\ka}}  \geq 
\delta  \right] \ ; \
\right. \nonumber \\ & & \mbox{~} \left.
\limsup_{n\tinf} \frac{v_n^2}{nh_n^{*d}} 
\bP  \left[  {\ti W}_n\cap [W_n(\zeta)]^c\neq\emptyset  \right]
\right\}.
\end{eqnarray*}
Moreover, following the proof of (\ref{a2refeq}) and 
(\ref{bis a2refeq}), we get
\begin{eqnarray*}
\limsup_{n\tinf}\frac{v_n^2}{nh_n^{*d}} 
\log \bP \, \left[ \,  {\ti U}_n\cap [U_n(\xi)]^c\neq\emptyset
\, \right]
& = &
-\infty ,
\nonumber
\\
\limsup_{n\tinf}
\frac{v_n^2}{nh_n^{*d}} 
\log \bP \, \left[ \,  {\ti W}_n\cap [W_n(\zeta)]^c\neq\emptyset
\, \right] 
& = &
-\infty .
\end{eqnarray*}
It thus follows from the application of Lemmas 
\ref{essai 1 lprobadusup} and \ref{essai 2 lprobadusup} 
that 
$$\limsup_{n\tinf} \frac{v_n^2}{nh_n^{*d}} 
\log \bP  \left(   \exists x\in C,\ \ f(x)\not\in
\left[f_n^*(x)-\delta\frac{\sqrt{\tilde T_n(x)\kappa}}{v_n}\ ;\
f_n^*(x)+\delta\frac{\sqrt{\tilde T_n(x)\kappa}}{v_n}\right]
 \right)
 \leq -\delta^2\left(1 -\frac{\zeta}{2}\right). $$
Since this last upper bound holds for any $\zeta>1$, 
it follows that
$$
\limsup_{n\tinf} \frac{v_n^2}{nh_n^{*d}} 
\log \bP  \left(   \exists x\in C,\ \ f(x)\not\in
\left[f_n^*(x)-\delta\frac{\sqrt{\tilde T_n(x)\kappa}}{v_n}\ ;\
f_n^*(x)+\delta\frac{\sqrt{\tilde T_n(x)\kappa}}{v_n}\right]
 \right)
 \leq -\frac{\delta^2}{2}. 
$$

\paragraph{Lower bound}

Following the proof of (\ref{ancien lminoration}), 
we obtain 
\begin{eqnarray*}
\liminf_{n\tinf} \frac{v_n^2}{nh_n^{*d}} 
\log \bP  \left(   \exists x\in C,\ \ f(x)\not\in
\left[f_n^*(x)-\delta\frac{\sqrt{\tilde T_n(x)\kappa}}{v_n}\ ;\
f_n^*(x)+\delta\frac{\sqrt{\tilde T_n(x)\kappa}}{v_n}\right]
 \right)
& \geq &
-\frac{\delta^2}{2},
\end{eqnarray*}
which concludes the proof of Theorem \ref{th2}.

\subsubsection{Proof of Theorem \ref{th1}}
\label{3.3.3}

\paragraph{Upper bound}

Let $(\e_n)$ be a sequence satisfying (\ref{assumption}) 
(the existence of
such a sequence is obvious in view of (\ref{b})). 
Moreover, set $C_n=C$ for all $n$, set $\xi\in]0,1[$, 
and let $U_n(\xi)$ be defined according to (\ref{4 global}).
We note that 
\begin{eqnarray*}
&&\bP \, \left( \,  \exists x\in C,\ f(x)\not\in
\left[f_n^*(x)-\delta\frac{\sqrt{f_n(x)\kappa}}{v_n}\ ;\
f_n^*(x)+\delta\frac{\sqrt{f_n(x)\kappa}}{v_n}\right]
\, \right)\\
&  & \leq
\bP \, \left[ \,  \sup_{x\in C}
\frac{v_n|f_n^*(x)-f(x)|}{\sqrt{f_n(x)\ka}} > \delta  
\ \ \mbox{and}\ \ 
\inf_{x\in C}f_n(x)>0\, \right]
+\bP\, \left[ \,  
\inf_{x\in C}f_n(x)=0\, \right]
\end{eqnarray*}
Since, under the assumptions of Theorem 
\ref{th1}, there exists $a>0$ such that 
$f(x)\geq a$ for all $x\in C$, we
have, for $n$ large enough, $U_n(\xi)=C$. It follows that  
\begin{eqnarray}
& &
\limsup_{n\tinf} \frac{v_n^2}{nh_n^{*d}} 
\log \bP  \left(   \exists x\in C,\ \ f(x)\not\in
\left[f_n^*(x)-\delta\frac{\sqrt{f_n(x)\kappa}}{v_n}\ ;\
f_n^*(x)+\delta\frac{\sqrt{f_n(x)\kappa}}{v_n}\right]
 \right) \nonumber \\
& &  \leq 
\max\left\{\limsup_{n\tinf} \frac{v_n^2}{nh_n^{*d}} 
\bP \, \left[  \sup_{x\in U_n(\xi)} 
 \frac{{v_n}\left|f_n^*(x)-f(x)\right|}{\sqrt{ f_n(x)\ka}}  
\geq   \delta 
\ \ \mbox{and}\ \ \inf_{x\in U_n(\xi)}f_n(x)>0 \right]\ ; \ 
\right. \nonumber \\ & & \mbox{~~} \left.
\limsup_{n\tinf} \frac{v_n^2}{nh_n^{*d}} 
\bP\, \left[ \,  
\inf_{x\in C}f_n(x)=0\, \right]
\right\},
\label{3Max}
\end{eqnarray}
with, by application of Corollary 2 
in Mokkadem, Pelletier and Worms (2005),
\begin{eqnarray*}
 \lefteqn{\lim_{n\tinf} \frac{v_n^2}{nh_n^{*d}}\log
\bP\, \left[ \,  \sup_{x\in C}|f_n(x)-f(x)|\geq a\, \right]}\\
& = &
\lim_{n\tinf} \frac{v_n^2h_n^d}{h_n^{*d}}\left[\frac{v_n^2}{nh_n^{*d}}\log
\bP\, \left[ \,  \sup_{x\in C}|f_n(x)-f(x)|\geq a\, \right]\right]\\
& = &
-\infty.
\end{eqnarray*}
The application of Lemma \ref{essai 1 lprobadusup} then gives 
$$
 \limsup_{n\tinf} \frac{v_n^2}{nh_n^{*d}} 
\log \bP \, \left( \,  \exists x\in C,\ f(x)\not\in
\left[f_n^*(x)-\delta\frac{\sqrt{f_n(x)\kappa}}{v_n}\ ;\
f_n^*(x)+\delta\frac{\sqrt{f_n(x)\kappa}}{v_n}\right]
\, \right)
\; \leq \;  -\frac{\delta^2}{2}.
$$

\paragraph{Lower bound}
Set $x_0\in C$; since $f(x_0)\neq 0$, we clearly have, 
by application of Theorem \ref{th0},
\begin{eqnarray*}
\lefteqn{\liminf_{n\tinf} \frac{v_n^2}{nh_n^{*d}} 
\log \bP \, \left( \,  \exists x\in C,\ f(x)\not\in
\left[f_n^*(x)-\delta\frac{\sqrt{f_n(x)\kappa}}{v_n}\ ;\
f_n^*(x)+\delta\frac{\sqrt{f_n(x)\kappa}}{v_n}\right]
\, \right) }\\
& \geq & 
\liminf_{n\tinf} \frac{v_n^2}{nh_n^{*d}} 
\log \bP \, \left( \,  f(x_0)\not\in
\left[f_n^*(x_0)-\delta\frac{\sqrt{f_n(x_0)\kappa}}{v_n}\ ;\
f_n^*(x_0)+\delta\frac{\sqrt{f_n(x_0)\kappa}}{v_n}\right]
\, \right) \\
& \geq & 
\frac{-\delta^2}{2},
\end{eqnarray*} 
which concludes the proof of Theorem \ref{th1}.

%
%
\subsection{Proof of Corollaries \ref{co1} and \ref{co2}}
\label{3.4}

Set $H=C$ in the framework of Corollary \ref{co1} and $H=\bR^d$ in the 
framework of Corollary \ref{co2}.
Let $\theta$ and $\theta_n$ satisfy 
$$f(\theta)=\sup_{x\in H}f(x)\ \ \mbox{and}\ \ 
f_n(\theta_n)=\sup_{x\in H}f_n(x),$$
and let $\ti T_n$ be the truncating function defined as 
$\ti T_n(x)=\max\left\{f_n(x);\e_nf(\theta)\right\}$ for all $x\in H$.

\paragraph{Upper bound}
We first prove that
\begin{equation}
\label{majocor}
 \limsup_{n\tinf} \frac{v_n^2}{nh_n^{*d}} 
\log \bP \, \left[ \, \exists{x\in H},\  
  {{v_n}|f_n^*(x)-f(x)|} \, > \,  
\delta{\sqrt{T_n(x)\kappa}} \, \right]
 \; \leq \;  -\frac{\delta^2}{2}.
\end{equation}
\\ 
Set $\eta>0$; since 
\begin{eqnarray*}
f_n(\theta)>\frac{f(\theta)}{1+\eta} & \imp & 
\forall x\in C,\ T_n(x)\geq \frac{\ti T_n(x)}{1+\eta} ,
\end{eqnarray*}
we have 
\begin{eqnarray*}
\lefteqn{\bP \, \left[ \, \exists{x\in H},\  
  {{v_n}|f_n^*(x)-f(x)|}\, > \,  
\delta{\sqrt{T_n(x)\kappa}}  \, \right]}\\
& \leq &
\bP \, \left[ \, \exists{x\in H},\  
  {{v_n}|f_n^*(x)-f(x)|}\, > \,  
\delta{\sqrt{T_n(x)\kappa}}
\ \ \mbox{and}\ \ f_n(\theta)>\frac{f(\theta)}{1+\eta}\, \right]
+
\bP\, \left[ \, f_n(\theta)\leq \frac{f(\theta)}{1+\eta}
\, \right]
\\
& \leq &
\bP \, \left[ \, \exists{x\in H},\  {{v_n}|f_n^*(x)-f(x)|}\, > \,  
\frac{\delta\sqrt{\ti T_n(x)\kappa}}{\sqrt{1+\eta}}\, \right]
+
\bP \, \left[ \, f(\theta)-f_n(\theta)\geq \frac{\eta f(\theta)}{1+\eta}
   \, \right].
\end{eqnarray*}
The application of Theorem \ref{th2} (respectively Theorem 
\ref{th3}) in the case $H=C$ (respectively $H=\bR^d$) ensures that
\begin{eqnarray*}
 \limsup_{n\tinf} \frac{v_n^2}{nh_n^{*d}} 
\log \bP \, \left[ \, \exists{x\in H},\  {{v_n}|f_n^*(x)-f(x)|}\, > \,  
\frac{\delta\sqrt{\ti T_n(x)\kappa}}{\sqrt{1+\eta}}\, 
 \right]
& \leq &  
-\frac{\delta^2}{2(1+\eta)}
\end{eqnarray*}
and the application of Corollary 1 in Mokkadem, Pelletier 
and Worms (2005) gives
\begin{eqnarray}
\lefteqn{\limsup_{n\tinf} \frac{v_n^2}{nh_n^{*d}} 
\log \bP \, \left[ \,  f(\theta)-f_n(\theta)\geq \frac{\eta f(\theta)}{1+\eta}
\, \right]} \nonumber\\
& \leq &
\limsup_{n\tinf} \left[\frac{v_n^2h_n^d}{nh_n^{*d}}\right]
\left[\frac{1}{nh_n^{d}} \log \bP \, \left[ \,  
f(\theta)-f_n(\theta)\geq \frac{\eta f(\theta)}{1+\eta}
\, \right]\right] \nonumber\\
& = & -\infty. \nonumber
\end{eqnarray}
Thus, we get  
\[
 \limsup_{n\tinf} \frac{v_n^2}{nh_n^{*d}} 
\log \bP \, \left[ \, \exists{x\in H},\  
  {{v_n}|f_n^*(x)-f(x)|} \, > \,  
\delta{\sqrt{T_n(x)\kappa}} \, \right]
 \; \leq \;  -\frac{\delta^2}{2(1+\eta)},
\]
and, since $\eta$ can be chosen arbitrarily close to zero, the proof 
of the upper bound (\ref{majocor}) is completed.
\paragraph{Lower bound} 
We now prove the lower bound
\begin{equation}
\label{minocor}
 \liminf_{n\tinf} \frac{v_n^2}{nh_n^{*d}} 
\log \bP \, \left[ \, \exists{x\in H},\  
  {{v_n}|f_n^*(x)-f(x)|} \, > \,  
\delta{\sqrt{T_n(x)\kappa}} \, \right]
 \; \geq \;  -\frac{\delta^2}{2}.
\end{equation}
Set $\eta>0$; since 
\begin{eqnarray*}
f_n(\theta_n)\leq (1+\eta){f(\theta)} & \imp & 
\forall x\in C,\ T_n(x)\leq (1+\eta){\ti T_n(x)} ,
\end{eqnarray*}
we have 
\begin{eqnarray*}
\lefteqn{\bP \, \left[ \, \exists{x\in H} ,\ 
  {{v_n}\left|f_n^*(x)-f(x)\right|}
  \, > \,  \delta{\sqrt{T_n(x)\ka}} \, \right]}\\
& \geq &
\bP \, \left[ \,  \exists{x\in H} ,\ 
  {{v_n}\left|f_n^*(x)-f(x)\right|}
  \, > \,  \delta{\sqrt{(1+\eta)\ti T_n(x)\ka}} 
\ \ \mbox{and}\ \ f_n(\theta_n)\leq (1+\eta){f(\theta)} \, \right]\\
& \geq &
\bP \, \left[ \, 
\exists{x\in H} ,\ 
  {{v_n}\left|f_n^*(x)-f(x)\right|}
  \, > \,  \delta{\sqrt{(1+\eta)\ti T_n(x)\ka}}  
 \, \right]
- \bP \, \left[ \,  f_n(\theta_n)-f(\theta)> \eta{f(\theta)}
   \, \right]\\
& \geq & 
\bP \, \left[ \, 
\exists{x\in H} ,\ 
  {{v_n}\left|f_n^*(x)-f(x)\right|}
  \, > \,  \delta{\sqrt{(1+\eta)\ti T_n(x)\ka}}  
 \, \right]
- \bP \, \left[ \,   \sup_{x\in H}|f_n(x)-f(x)|> \eta{f(\theta)}
   \, \right].
\end{eqnarray*}
The application of Theorem \ref{th2} (respectively Theorem 
\ref{th3}) in the case $H=C$ (respectively $H=\bR^d$) ensures that
\begin{eqnarray*}
 \liminf_{n\tinf} \frac{v_n^2}{nh_n^{*d}} 
\log \bP \, \left[ \, \exists{x\in H},\  {{v_n}|f_n^*(x)-f(x)|}\, > \,  
\delta{\sqrt{(1+\eta)\ti T_n(x)\ka}}  \, 
 \right]
& \geq &  
-\frac{\delta^2(1+\eta)}{2}
\end{eqnarray*}
and the application of Corollary 2 in Mokkadem, Pelletier 
and Worms (2005) ensures that
\begin{eqnarray*}
\lefteqn{
\limsup_{n\tinf} \frac{v_n^2}{nh_n^{*d}} 
\log \bP \, \left[ \,  \sup_{x\in H}|f_n(x)-f(x)|> \eta{f(\theta)}  \, \right]
}\\
& = &
\limsup_{n\tinf} \frac{v_n^2h_n^d}{h_n^{*d}} 
\left\{ \frac{1}{nh_n^{d}} 
\log \bP \, \left[ \,   \sup_{x\in H}|f_n(x)-f(x)|> \eta{f(\theta)}
\, \right]\right\}\\
& = &
-\infty.
\end{eqnarray*}
Thus, it follows that 
\begin{eqnarray*}
 \liminf_{n\tinf} \frac{v_n^2}{nh_n^{*d}} 
\log \bP \, \left[ \, \exists{x\in H},\  
  {{v_n}|f_n^*(x)-f(x)|} \, > \,  
\delta{\sqrt{T_n(x)\kappa}} \, \right]
 & \geq &  -\frac{\delta^2}{2}, 
\end{eqnarray*}
and, since $\eta$ can be chosen arbitrarily close to zero, 
the lower bound (\ref{minocor}) follows.

\section{Appendix}

\subsection{Proof of Lemmas \ref{essai 1 lprobadusup} 
and \ref{essai 2 lprobadusup}}
\label{Sec App 1}

The proof of Lemmas \ref{essai 1 lprobadusup} 
and \ref{essai 2 lprobadusup} requires the two following
preliminary lemmas.

\begin{lemma}
\label{lxisupedesproba}
Under Assumptions (A1)-(A3), we have
$$
 \limsup_{n\tinf} \frac{v_n^2}{nh_n^{*d}} \sup_{x\in U_n(\xi)}
\log \bP \, \left[ \,  \frac{v_n|f_n^*(x)-f(x)|}{\sqrt{f(x)\ka}}
 \geq \delta  \, \right] \; \leq \;  -\frac{\delta^2}{2}.
$$
\end{lemma}
\vspace{0.5cm}

\begin{lemma}
\label{lzetasupedesproba}
Under Assumptions (A1)-(A3), we have
$$
 \limsup_{n\tinf} \frac{v_n^2}{nh_n^{*d}} \sup_{x\in W_n(\zeta)}
\log \bP \, \left[ \,  \frac{v_n|f_n^*(x)-f(x)|}{\sqrt{\e_n\kappa}}
 \geq \delta  \, \right] \; \leq \; -\delta^2\left(1-\frac{\zeta}{2}\right).
$$
\end{lemma}
\vspace{0.5cm}

We first prove Lemmas \ref{lxisupedesproba} and 
\ref{lzetasupedesproba} in Subsection 
\ref{Subsection 3.5.1},
and then establish Lemmas \ref{essai 1 lprobadusup} 
and \ref{essai 2 lprobadusup} in Subsection 
\ref{Subsection 3.5.2}.

\subsubsection{Proof of Lemmas \ref{lxisupedesproba} and 
\ref{lzetasupedesproba}}

\label{Subsection 3.5.1}

Set 
\begin{eqnarray*}
u_n(x) & = & 
\left\{\begin{array}{ll}
\frac{\de}{\sqrt{f(x)\kappa}} & 
\mbox{in the framework of Lemma \ref{lxisupedesproba} }  \\
\frac{\de}{\sqrt{\e_n\kappa}}
& \mbox{in the framework of Lemma \ref{lzetasupedesproba}} 
\end{array}\right. \\
{\cal E} & = & 
\left\{\begin{array}{ll}
U_n(\xi) & \mbox{in the framework of Lemma \ref{lxisupedesproba} }  \\
W_n(\zeta) & \mbox{in the framework of Lemma \ref{lzetasupedesproba}} 
\end{array}\right.
\end{eqnarray*}
and, for any $u\in\bR$,
\begin{eqnarray*}
 \Lambda_{n,x}(u) & = & \frac{v_n^2}{nh_n^{*d}} 
\log \bE \, \left[ \, \exp\left(\frac{nh_n^{*d}}{v_n}u
[f_n^*(x)-f(x)]\right) \,\right].    
\end{eqnarray*}
\\
To study the asymptotics of 
$\sup_{x\in {\cal E}}\bP \, \left[ \,  v_n u_n(x)|f_n^*(x)-f(x)| 
\geq \delta^2  \, \right]$,  we first note that, by Chebyshev's inequality,
we have  
\begin{eqnarray*}
\lefteqn{\bP \, [ \, v_n u_n(x)[f_n^*(x)-f(x)] \geq \delta^2 \, ]  }\\
& = &
\bP \, \left[ \, \exp\left(\frac{nh_n^{*d}}{v_n} u_n(x)[f_n^*(x)-f(x)]\right) 
\geq \exp\left(\frac{nh_n^{*d}}{v_n^2}\delta^2\right) \, \right]  \\ 
& \leq &
\exp \left[-\frac{nh_n^{*d}}{v_n^2}\delta^2\right]
\bE\left(\exp \left[\frac{nh_n^{*d}}{v_n}u_n(x)[f_n^*(x)-f(x)]
\right]\right)\\  
& \leq &
\exp \left[-\frac{nh_n^{*d}}{v_n^2}\delta^2\right]
\exp \left[\frac{nh_n^{*d}}{v_n^2}\Lambda_{n,x}(u_n(x))\right]
\end{eqnarray*}
and thus
\begin{equation}
\label{aaa}
 \frac{v_n^2}{nh_n^{*d}} 
\log \sup_{x\in {\cal E}} 
\bP \, [ \,v_n u_n(x)[f_n^*(x)-f(x)] \geq \delta^2 \, ] \; \leq \;
 - \delta^2 + \sup_{x\in{\cal E}} \Lambda_{n,x}(u_n(x)).
\end{equation}
In the same way, we prove that
\begin{equation}
\label{bbb}
 \frac{v_n^2}{nh_n^{*d}} 
\log \sup_{x\in {\cal E}} 
\bP \, [ \,v_n u_n(x)[f(x)-f_n^*(x)] \geq \delta^2 \, ] \; \leq \;
 - \delta^2 + \sup_{x\in{\cal E}} \Lambda_{n,x}(-u_n(x)).
\end{equation}
\\
Let us set $e\in\{-1,+1\}$ and let us at first assume that
\begin{equation}
\label{cul}
 \Lambda_{n,x}(eu_n(x))  =
\frac{u_n^2(x)}{2}\ka f(x)+R_{n,x}(eu_n(x)) \ \ \mbox{with}\ \ 
\lim_{n\tinf}\sup_{x\in{\cal E}}R_{n,x}(eu_n(x))=0.
\end{equation}
\\
$\bullet$
In the framework of Lemma \ref{lxisupedesproba}, (\ref{cul}) means that 
$$ \Lambda_{n,x}(eu_n(x))  =
\frac{\de^2}{2}+R_{n,x}(eu_n(x)) \ \ \mbox{with}\ \ 
\lim_{n\tinf}\sup_{x\in U_n(\xi)}R_{n,x}(eu_n(x))=0, $$
and Lemma \ref{lxisupedesproba} is thus a straightforward consequence of
(\ref{aaa}) and (\ref{bbb}).\\
$\bullet$
In the framework of Lemma \ref{lzetasupedesproba}, (\ref{cul}) can be rewritten
as 
$$ \Lambda_{n,x}(eu_n(x))  =
\frac{\de^2}{2\e_n}f(x)+R_{n,x}(eu_n(x)) \ \ \mbox{with}\ \ 
\lim_{n\tinf}\sup_{x\in W_n(\zeta)}R_{n,x}(eu_n(x))=0,  $$
and, since $\sup_{x\in W_n(\zeta)}f(x)/\e_n\leq\zeta$, Lemma 
\ref{lzetasupedesproba} is also given by (\ref{aaa}) and (\ref{bbb}).
\\
\\
It remains to prove (\ref{cul}). Let us first note that
\begin{eqnarray*}
 \Lambda_{n,x}(eu_n(x)) & = &
-v_neu_n(x)f(x)+  
\frac{v_n^2}{nh_n^{*d}}\log\bE \, \left[ \,
  \exp \left(  \, v_n^{-1}eu_n(x)
\sum_{i=1}^nK\left(\frac{x-X_i}{h_n^*}\right)
\, \right) \, \right]\\
 & = &
-v_neu_n(x)f(x)+
\frac{v_n^2}{h_n^{*d}}  \log  \bE \, \left[ \, 
\exp\left\{{v_n^{-1}eu_n(x)K\left(\frac{x-X_1}{h_n^*}\right)}\right\} 
\, \right] 
\end{eqnarray*}
By Taylor formula for the function $\log$, there exists $c_n$ between $1$ and 
$\bE\left[ \, 
\exp\left\{{v_n^{-1}eu_n(x)K\left(\frac{x-X_1}{h_n^*}\right)}\right\} 
\, \right]$ such that
\begin{eqnarray*}
 \Lambda_{n,x}(eu_n(x))
 & = &
-v_neu_n(x)f(x)+
\frac{v_n^2}{h_n^{*d}}  \bE \, \left[ \, 
\exp\left\{{v_n^{-1}eu_n(x)K\left(\frac{x-X_1}{h_n^*}\right)}\right\} 
-1\, \right]
 - R^{(1)}_{n,x}(eu_n(x))
\end{eqnarray*}
with
\begin{eqnarray}
\label{resteun}
 R^{(1)}_{n,x}(eu_n(x))  & = &
 \frac{v_n^2}{2c_n^2h_n^{*d}} 
{\left\{\bE \, \left[ \, 
\exp\left\{{v_n^{-1}eu_n(x)K\left(\frac{x-X_1}{h_n^*}\right)}\right\} 
-1\, \right]
\right\}}^2.
\end{eqnarray}
We rewrite $\Lambda_{n,x}(eu_n(x))$ as
\begin{eqnarray}
\lefteqn{\Lambda_{n,x}(eu_n(x))}\nonumber \\ 
& = &
-v_neu_n(x)f(x)+
\frac{v_n^2}{h_n^{*d}} \int_{\bR^d} \left\{
 \exp\left[v_n^{-1}eu_n(x)K\left(\frac{x-y}{h^*_n}\right)\right] -1
  \right\} f(y) \, dy  -R^{(1)}_{n,x}(eu_n(x)) \nonumber  \\
 & = &
-v_neu_n(x)f(x)+
 v_n^2 \int_{\bR^d} \left(\exp\left[v_n^{-1}eu_n(x)K(z)\right] -1 \right) 
f(x-h_n^*z) \, dz
 -R^{(1)}_{n,x}(eu_n(x))
 \nonumber\\
 & = &
-v_neu_n(x)f(x)+ 
v_n^2 \int_{\bR^d} \left(  
v_n^{-1}eu_n(x)K(z)+\frac{v_n^{-2}u_n^2(x)K^2(z)}{2}
\right) f(x-h^*_nz) \, dz \nonumber\\
 &  & \mbox{}
 -R^{(1)}_{n,x}(eu_n(x)) + R^{(2)}_{n,x}(eu_n(x))
\nonumber
\end{eqnarray}
with
\begin{eqnarray}
 R^{(2)}_{n,x}(eu_n(x))
 & \leq &
 v_n^2 \int_{\bR^d}\frac{v_n^{-3}u_n^3(x)|K(z)|^3}{6} 
\exp\left[v_n^{-1}u_n(x)|K(z)|\right]
f(x-h_n^*z) \, dz \nonumber\\
 & \leq &
\frac{v_n^{-1}u_n^3(x)}{6}\|f\|_{\infty}
\exp\left[v_n^{-1}u_n(x)\|K\|_{\infty}\right]\int_{\bR^d}|K(z)|^3\, dz.
\label{restedeux}
\end{eqnarray}
It follows that
\begin{eqnarray}
\lefteqn{\Lambda_{n,x}(eu_n(x))}\nonumber \\ 
& = &
-v_neu_n(x)\int_{\bR^d} K(z)\left[f(x)-f(x-h^*_nz)\right] \, dz  
+\frac{u_n^2(x)}{2}  \int_{\bR^d} K^2(z)f(x-h_n^*z) \, dz \nonumber\\
 &  & \mbox{}
 -R^{(1)}_{n,x}(eu_n(x)) + R^{(2)}_{n,x}(eu_n(x))\nonumber  
\end{eqnarray}
and, setting 
\begin{equation}
\label{restetrois}
R^{(3)}_{n,x}(eu_n(x))=v_neu_n(x)
\int_{\bR^d} K(z)\left[f(x)-f(x-h^*_nz)\right] \, dz,
\end{equation}
we obtain
\begin{eqnarray}
\lefteqn{\Lambda_{n,x}(eu_n(x))}\nonumber \\ 
& = &
\frac{u_n^2(x)}{2}  \int_{\bR^d} K^2(z)\left[f(x-h_n^*z)-f(x)\right] \, dz
+\frac{u_n^2(x)f(x)\ka}{2} -R^{(1)}_{n,x}(eu_n(x)) 
\nonumber\\ &  & \mbox{}
 + R^{(2)}_{n,x}(eu_n(x)) -R^{(3)}_{n,x}(eu_n(x))\nonumber \\ 
& = &
\frac{u_n^2(x)f(x)\ka}{2}-R^{(1)}_{n,x}(eu_n(x)) 
+ R^{(2)}_{n,x}(eu_n(x)) -R^{(3)}_{n,x}(eu_n(x))
+ R^{(4)}_{n,x}(eu_n(x)) \nonumber
\end{eqnarray}
where 
\begin{equation}
\label{restequatre}
R^{(4)}_{n,x}(eu_n(x))=\frac{u_n^2(x)}{2}
\int_{\bR^d} K^2(z)\left[f(x-h_n^*z)-f(x)\right] \, dz.
\end{equation}
To conclude the proof of Lemmas \ref{lxisupedesproba} and 
\ref{lzetasupedesproba}, it remains to show that
$$\lim_{n\tinf}\sup_{x\in{\cal E}}R^{(i)}_{n,x}(eu_n(x))=0\ \ \mbox{for}\ \ 
i\in\{1,\ldots ,4\}.$$
Let $c$ and $c'$ denote generic positive constants that
may vary from line to line. We shall use several times the fact that
$$\sup_{x\in{\cal E}}u_n(x)\leq\frac{c}{\sqrt{\e_n}}.$$
\\
%
%
%
$\bullet$
To study $R^{(1)}_{n,x}(eu_n(x))$ defined in (\ref{resteun}), we first 
note that 
$$\exp\left\{{v_n^{-1}eu_n(x)K\left(\frac{x-X_1}{h_n^*}\right)}\right\} 
\geq\exp\left\{-v_n^{-1}u_n(x)\|K\|_{\infty}\right\} $$
so that
$$
 \frac{1}{c_n^2}\leq\exp\left\{2v_n^{-1}u_n(x)\|K\|_{\infty}\right\}.
$$
Moreover, since
\begin{eqnarray}
\lefteqn{ \left| \bE \,  
\left[ \, 
\exp\left\{{v_n^{-1}eu_n(x)K\left(\frac{x-X_1}{h_n^*}\right)}\right\} 
-1\, \right]
\,  \right|} \nonumber \\
 & \leq &
 \bE \, \left[ \, 
\left|v_n^{-1}eu_n(x)K\left(\frac{x-X_1}{h_n^*}\right)\right|
\exp\left|{v_n^{-1}eu_n(x)K\left(\frac{x-X_1}{h_n^*}\right)}\right|
 \, \right] \nonumber \\
 & \leq &
 h_n^{*d}v_n^{-1}u_n(x) \|f\|_{\infty} 
\exp\left\{v_n^{-1}u_n(x)\|K\|_{\infty}\right\}
{\tts\int_{\bR^d}} |K(z)|\, dz,
 \nonumber
\end{eqnarray}
we obtain
\begin{eqnarray*}
\sup_{x\in {\cal E}}
 |R^{(1)}_{n,x}(eu_n(x))| 
& \leq &
\sup_{x\in {\cal E}}  h_n^{*d} \frac{u_n^2(x)}{2} 
\exp\left\{4v_n^{-1}u_n(x)\|K\|_{\infty}\right\}
\|f\|_{\infty}^2
 \left( {\tts\int_{\bR^d}}|K(z)|\, dz\right)^2\\
& \leq &
c\frac{h_n^{*d}}{\e_n}\exp\left(\frac{c'}{v_n\sqrt{\e_n}}\right) \\
& \ra & 0\ \ \mbox{since} \ \ 
h_n^{*d}/\e_n\ra 0\ \ \mbox{and} \ \ v_n\sqrt{\e_n}\tinf .
\end{eqnarray*}
\\
%
$\bullet$
It follows from (\ref{restedeux}) that
\begin{eqnarray*}
\sup_{x\in {\cal E}} |R^{(2)}_{n,x}(eu_n(x))|
& \leq &
\frac{c}{v_n\e_n^{3/2}}\exp\left(\frac{c}{v_n\sqrt{\e_n}}\right) \\
& \ra & 0\ \ \mbox{since} \ \  v_n\e_n^{3/2}\tinf .
\end{eqnarray*}
\\
%
$\bullet$
To upper bound $R^{(3)}_{n,x}$ defined in (\ref{restetrois}), 
we use a Taylor expansion and Assumptions (A1)-(A3) to obtain
$$\sup_{x\in {\cal E}}\left|
\int_{\bR^d} K(z)\left[f(x)-f(x-h^*_nz)\right] \, dz\right|
\leq h_n^{*2}\sup_{x\in\bR^d}\|D^2f(x)\|\int_{\bR^d}\|z\|^2K(z)dz,$$
from which we deduce that
\begin{eqnarray*}
\sup_{x\in {\cal E}} |R^{(3)}_{n,x}(eu_n(x))|
& \leq &
c\frac{v_nh_n^{*2}}{\sqrt{\e_n}} \\
& \ra & 0\ \ \mbox{since} \ \ 
\frac{v_nh_n^{*2}}{\sqrt{\e_n}}\ra 0 .
\end{eqnarray*}
\\
%
$\bullet$
Similarly, for $R^{(4)}_{n,x}$ defined in (\ref{restequatre}), 
we note that 
$$\sup_{x\in {\cal E}}\left|
\int_{\bR^d} K^2(z)\left[f(x)-f(x-h^*_nz)\right] \, dz\right|
\leq h_n^{*}\|K\|_\infty 
\sup_{x\in\bR^d}\|\nabla f(x)\|\int_{\bR^d}\|z\|K(z)dz,$$
so that
\begin{eqnarray*}
\sup_{x\in {\cal E}} |R^{(4)}_{n,x}(eu_n(x))|
& \leq &
c\frac{h_n^{*}}{\e_n} \\
& \ra & 0\ \ \mbox{since} \ \ 
\frac{h_n^{*}}{{\e_n}}\ra 0 , 
\end{eqnarray*}
which concludes the proof of Lemmas \ref{lxisupedesproba} and 
\ref{lzetasupedesproba}.

\subsubsection{Proof of Lemmas \ref{essai 1 lprobadusup} 
and \ref{essai 2 lprobadusup}}
\label{Subsection 3.5.2}

In view of Assumption (A2), there exists $\beta$ and $\|K\|_H$ such that
$$|K(x)-K(y)|\leq\|K\|_H\|x-y\|^\beta\ \ \forall x,y\in\bR^d.$$
Set $b=\sup_{x\in\bR^d}\|\nabla f(x)\|$, $\delta>0$, 
$\rho \in ] 0,\delta[$, and
\begin{eqnarray*}
 R_n & = &
\left( \frac{\rho h_n^{*(d+\beta)}\sqrt{\xi\e_n\kappa}}
 {2[2\sqrt{d}]^\beta v_n(b+\|K\|_H)} \right)^{1/\beta}.
\end{eqnarray*}
Let $m$ be the integer satisfying $w_n/R_n\leq m<1+w_n/R_n$, and set 
$N'(n)=m^d$. The ball $D_n=\left\{x\in\bR^d,\ \|x\|\leq w_n\right\}$ can
be covered by $N'(n)$ cubes with lenght side $2R_n$ (note that $ U_n(\xi)$,
$W_n(\zeta)$, and $C_n$ are subsets of $D_n$). We denote by
$B^{(i)}_n$, $i=1,\ldots ,N(n)$ ($N(n)\leq N'(n)$), the cubes that intersect
$U_n(\xi)$, and by $\ti B^{(i)}_n$, $i=1,\ldots ,\ti N(n)$ 
($\ti N(n)\leq N'(n)$), the cubes that intersect $W_n(\zeta)$. 
Moreover, for each $i\in\{1,\ldots ,N(n)\}$, we choose 
$x^{(i)}_n\in B^{(i)}_n\cap U_n(\xi)$, and for each 
$i\in\{1,\ldots ,\ti N(n)\}$, we choose 
$\ti x^{(i)}_n\in \ti B^{(i)}_n\cap W_n(\zeta)$.\\ 
\\
To prove Lemma \ref{essai 1 lprobadusup}, we first 
note that
\begin{eqnarray*}
 \bP \, \left[ \, \sup_{x\in U_n(\xi)}
\frac{{v_n}\left|f_n^*(x)-f(x)\right|}{\sqrt{f(x)\ka}}
\geq \delta \, \right]
 & \leq &
 \sum_{i=1}^{N(n)} \bP \, \left[ \, \sup_{x\in B^{(n)}_i\cap U_n(\xi)}
\frac{{v_n}\left|f_n^*(x)-f(x)\right|}{\sqrt{f(x)\ka}} 
\geq \delta  \, \right]  \\
 & \leq &
N(n)\max_{1\leq i\leq N(n)}\bP \, \left[ \, 
\sup_{x\in B^{(n)}_i\cap U_n(\xi)}
\frac{{v_n}\left|f_n^*(x)-f(x)\right|}{\sqrt{f(x)\ka}} 
\geq \delta  \, \right].
\end{eqnarray*}
Now, for any $i\in\{1,\ldots ,N(n)\}$ and for any 
$x\in B^{(n)}_i\cap U_n(\xi)$, we write
\begin{eqnarray*}
\frac{{v_n}\left|f_n^*(x)-f(x)\right|}{\sqrt{f(x)\ka}}
& \leq &
\frac{{v_n}\left|f_n^*(x)-f_n^*(x^{(n)}_i)\right|}{\sqrt{f(x)\ka}}
+\frac{{v_n}\left|f_n^*(x^{(n)}_i)-f(x^{(n)}_i)\right|}{\sqrt{f(x)\ka}}
+\frac{{v_n}\left|f(x^{(n)}_i)-f(x)\right|}{\sqrt{f(x)\ka}}.
\end{eqnarray*}
On the one hand, since $K$ is 
H\"older continuous and since $x\in B^{(n)}_i\imp
\|x-x^{(n)}_i\|\leq 2\sqrt{d}R_n$, we get
\begin{eqnarray*}
\frac{{v_n}\left|f_n^*(x)-f_n^*(x^{(n)}_i)\right|}{\sqrt{f(x)\ka}}
& \leq &
\frac{{v_n}\|K\|_H}{\sqrt{\xi\e_n\ka}\ h_n^{*d}}
{\left\|\frac{x-x^{(n)}_i}{h_n^*}\right\|}^\beta\\
& \leq &
\frac{{v_n}\|K\|_H[2\sqrt{d}R_n]^\beta}
{h_n^{*(d+\beta)}\sqrt{\xi\e_n\ka}} \\
& \leq &
\frac{\rho}{2}.
\end{eqnarray*}
On the other hand, we have
\begin{eqnarray*} 
\frac{{v_n}\left|f(x^{(n)}_i)-f(x)\right|}{\sqrt{f(x)\ka}}
& \leq &
\frac{{v_n}b
\left\|x-x^{(n)}_i\right\|}{\sqrt{\xi\e_n\ka}},
\end{eqnarray*}
and, since $\beta\leq 1$ and $R_n\ra 0$, we obtain, for $n$ large enough,
\begin{eqnarray*} 
\frac{{v_n}\left|f(x^{(n)}_i)-f(x)\right|}{\sqrt{f(x)\ka}}
& \leq &
\frac{{v_n}b
\left\|x-x^{(n)}_i\right\|}{\sqrt{\xi\e_n\ka}}\\
& \leq &
\frac{{v_n}bR_n^\beta}{\sqrt{\xi\e_n\ka}} \\
& \leq &
\frac{\rho}{2}. 
\end{eqnarray*}
We deduce that, for all $n$ sufficiently large,
$\forall i\in\{1,\ldots ,N(n)\}$,
$\forall x\in B^{(n)}_i\cap U_n(\xi)$,  
\begin{eqnarray*} 
\frac{{v_n}\left|f_n^*(x)-f(x)\right|}{\sqrt{f(x)\ka}}
& \leq &
\frac{{v_n}\left|f_n^*(x^{(n)}_i)-f(x^{(n)}_i)\right|}{\sqrt{f(x)\ka}}
+\rho\\
& \leq &
\sqrt{1+\frac{f(x^{(n)}_i)-f(x)}{f(x)}}\ 
\frac{{v_n}\left|f_n^*(x^{(n)}_i)-f(x^{(n)}_i)\right|}{\sqrt{
f(x^{(n)}_i)\ka}}+\rho\\
& \leq &
\sqrt{1+\frac{bR_n}{f(x)}}\ 
\frac{{v_n}\left|f_n^*(x^{(n)}_i)-f(x^{(n)}_i)\right|}{\sqrt{
f(x^{(n)}_i)\ka}}+\rho\\
& \leq &
\sqrt{1+\frac{bR_n}{\xi\e_n}}\ 
\frac{{v_n}\left|f_n^*(x^{(n)}_i)-f(x^{(n)}_i)\right|}{\sqrt{
f(x^{(n)}_i)\ka}}+\rho \\
& \leq &
\sqrt{1+\rho}\ 
\frac{{v_n}\left|f_n^*(x^{(n)}_i)-f(x^{(n)}_i)\right|}{\sqrt{
f(x^{(n)}_i)\ka}}+\rho.
\end{eqnarray*}
We can then deduce that, for $n$ large enough,
\begin{eqnarray*}
 \bP \, \left[ \, \sup_{x\in  U_n(\xi)}
\frac{{v_n}\left|f_n^*(x)-f(x)\right|}{\sqrt{f(x)\ka}}
\geq \delta \, \right]
 & \leq &
N(n)\max_{1\leq i\leq N(n)}\bP \, \left[ \, 
\frac{{v_n}\left|f_n^*(x^{(n)}_i)-f(x^{(n)}_i)\right|}
{\sqrt{f(x^{(n)}_i)\ka}} 
\geq \frac{\delta-\rho}{\sqrt{1+\rho}}  \, \right]\\
 & \leq &
N(n)\sup_{x\in  U_n(\xi)}\bP \, \left[ \, 
\frac{{v_n}\left|f_n^*(x)-f(x)\right|}
{\sqrt{f(x)\ka}} 
\geq \frac{\delta-\rho}{\sqrt{1+\rho}}  \, \right].
\end{eqnarray*}
Applying Lemma \ref{lxisupedesproba}, we obtain
\begin{eqnarray*}
\lefteqn{\limsup_{n\tinf} \frac{v_n^2}{nh_n^{*d}} 
\log \bP \, \left[ \,  \sup_{x\in U_n(\xi)}
\frac{v_n|f_n^*(x)-f(x)|}{\sqrt{f(x)\kappa}}
 \geq \delta  \, \right] }\\
& \leq &
\limsup_{n\tinf}\left\{\frac{v_n^2\log N(n)}{nh_n^{*d}} \right\}
-\frac{(\de-\rho)^2}{2(1+\rho)}\\
& \leq &
\limsup_{n\tinf}\left\{
\frac{dv_n^2}{\beta nh_n^{*d}}\left[\beta\log (w_n)  
-(d+\beta)\log h_n^*-\frac{\log\e_n}{2}+
\log v_n\right] \right\}
-\frac{(\de-\rho)^2}{2(1+\rho)}\\
& \leq &
-\frac{(\de-\rho)^2}{2(1+\rho)}
\ \ \left(\mbox{since}\ \ \frac{v_n^2\log h_n^*}{nh_n^{*d}}\ra 0,\ \ 
\frac{v_n^2\log v_n}{nh_n^{*d}}\ra 0,\ \ 
\ \ \mbox{and}\ \ \frac{v_n^2\log \e_n}{nh_n^{*d}}\ra 0\right).
\end{eqnarray*}
This last upper bound holding for any $\rho\in]0,\delta[$, it follows that
\begin{eqnarray}
\label{avecf}
\limsup_{n\tinf} \frac{v_n^2}{nh_n^{*d}} 
\log \bP \, \left[ \,  \sup_{x\in U_n(\xi)}
\frac{v_n|f_n^*(x)-f(x)|}{\sqrt{f(x)\kappa}}
 \geq \delta  \, \right] 
& \leq &
-\frac{\de^2}{2}.
\end{eqnarray}
\\
To conclude the proof of Lemma \ref{essai 1 lprobadusup}, 
let us now set $\eta>0$, and note that
\begin{eqnarray}
\lefteqn{\bP \, \left[ \, \sup_{x\in U_n(\xi)} 
  \frac{{v_n}\left|f_n^*(x)-f(x)\right|}{\sqrt{ f_n(x)\ka}}  \, \geq \,  
\delta 
\ \ \mbox{and}\ \ \inf_{x\in U_n(\xi)}f_n(x)>0
\, \right]} \nonumber\\
& \leq &
\bP \, \left[ \, \sup_{x\in U_n(\xi)} 
\frac{f(x)}{f_n(x)}  \, \geq \,  1+\eta 
\ \ \mbox{and}\ \ \inf_{x\in U_n(\xi)}f_n(x)>0
\, \right]\nonumber\\
& & \mbox{~}
+
\bP \, \left[ \, \sup_{x\in U_n(\xi)} 
  \frac{{v_n}\left|f_n^*(x)-f(x)\right|}{\sqrt{ f(x)\ka}}  \, \geq \,  
\frac{\delta}{\sqrt{1+\eta}} 
\, \right] \nonumber\\
& \leq &
\bP \, \left[ \, \sup_{x\in U_n(\xi)} 
\frac{f(x)-f_n(x)}{f(x)}  \, \geq \,  \frac{\eta}{1+\eta} \, \right]
+
\bP \, \left[ \, \sup_{x\in  U_n(\xi)} 
  \frac{{v_n}\left|f_n^*(x)-f(x)\right|}{\sqrt{ f(x)\ka}}  \, \geq \,  
\frac{\delta}{\sqrt{1+\eta}} \, \right].
\nonumber
\end{eqnarray}
Since $x\in  U_n(\xi)\imp f(x)\geq \xi\e_n$ and since 
$v_n^2h_n^d\e_n^2/h_n^{*d}\tinf$, we have, by application of 
Theorem 5 in Mokkadem, Pelletier and Worms (2005),
\begin{eqnarray*}
\lefteqn{\limsup_{n\tinf}\frac{v_n^2}{nh_n^{*d}}
\bP \, \left[ \, \sup_{x\in U_n(\xi)} 
\frac{f(x)-f_n(x)}{f(x)}  \, \geq \,  \frac{\eta}{1+\eta} \, \right]}\\
& \leq &
\limsup_{n\tinf}
\frac{v_n^2h_n^{d}\e_n^2}{h_n^{*d}}
\left\{\frac{1}{nh_n^{d}\e_n^2}\bP \, \left[ \, \frac{1}{\e_n}\sup_{x\in C_n} 
\left|f_n(x)-f(x)\right|  \, \geq \,  
\frac{\xi\eta}{1+\eta} \, \right]\right\}\\
& = &
-\infty. 
\end{eqnarray*}
Now, in view of (\ref{avecf}), we deduce that
\begin{eqnarray*}
\limsup_{n\tinf}\frac{v_n^2}{nh_n^{*d}} \log 
\bP \, \left[ \, \sup_{x\in U_n(\xi)} 
  \frac{{v_n}\left|f_n^*(x)-f(x)\right|}{\sqrt{ f_n(x)\ka}}  \, \geq \,  
\delta \, \right]
& \leq &
\frac{-\de^2}{2(1+\eta)},
\end{eqnarray*}
and, since this last upper bound holds for all $\eta>0$, 
Lemma \ref{essai 1 lprobadusup} follows.
\\ 
\\
To prove Lemma \ref{essai 2 lprobadusup}, we proceed 
exactly in the same way as for establihing (\ref{avecf});
we first note that 
\begin{eqnarray*}
 \bP \, \left[ \, \sup_{x\in W_n(\zeta)}
\frac{{v_n}\left|f_n^*(x)-f(x)\right|}{\sqrt{\e_n\ka}}
\geq \delta \, \right]
 & \leq &
\ti N(n)\max_{1\leq i\leq \ti N(n)}\bP \, \left[ \, \sup_{x\in 
\ti B^{(n)}_i\cap W_n(\zeta)}
\frac{{v_n}\left|f_n^*(x)-f(x)\right|}{\sqrt{\e_n\ka}} 
\geq \delta  \, \right].
\end{eqnarray*}
with, for any $i\in\{1,\ldots ,\ti N(n)\}$,  
$x\in \ti B^{(n)}_i\cap W_n(\zeta)$, and $n$ large enough,
\begin{eqnarray*}
\frac{{v_n}\left|f_n^*(x)-f(x)\right|}{\sqrt{\e_n\ka}}
& \leq &
\frac{{v_n}\left|f_n^*(x)-f_n^*(\ti x^{(n)}_i)\right|}{\sqrt{\e_n\ka}}
+\frac{{v_n}\left|f_n^*(\ti x^{(n)}_i)-f(\ti x^{(n)}_i)\right|}
{\sqrt{\e_n\ka}}
+\frac{{v_n}\left|f(\ti x^{(n)}_i)-f(x)\right|}{\sqrt{\e_n\ka}}\\
& \leq &
\frac{{v_n}\left|f_n^*(x^{(n)}_i)-f(x^{(n)}_i)\right|}{\sqrt{\e_n\ka}}
+\rho
\end{eqnarray*}
We then deduce that
\begin{eqnarray*}
 \bP \, \left[ \, \sup_{x\in W_n(\zeta)}
\frac{{v_n}\left|f_n^*(x)-f(x)\right|}{\sqrt{\e_n\ka}}
\geq \delta \, \right]
 & \leq &
\ti N(n)\sup_{x\in W_n(\zeta)}\bP \, \left[ \, 
\frac{{v_n}\left|f_n^*(x)-f(x)\right|}
{\sqrt{\e_n\ka}} 
\geq \delta-\rho  \, \right]
\end{eqnarray*}
and, applying  Lemma \ref{lzetasupedesproba}, we obtain
\begin{eqnarray*} 
\lefteqn{\limsup_{n\tinf} \frac{v_n^2}{nh_n^{*d}} 
\log \bP \, \left[ \,  \sup_{x\in W_n(\zeta)}
\frac{v_n|f_n^*(x)-f(x)|}{\sqrt{\e_n\kappa}}
 \geq \delta  \, \right] }\\
&\leq & 
\limsup_{n\tinf}\frac{v_n^2\log \ti N(n)}{nh_n^{*d}} 
-(\de-\rho)^2\left(1-\frac{\zeta}{2}\right)\\
&\leq & 
-(\de-\rho)^2\left(1-\frac{\zeta}{2}\right).
\end{eqnarray*}
Since this last upper bound holds for any $\rho>0$, 
Lemma \ref{essai 2 lprobadusup} follows.

\subsection{Proof of (\ref{ancien deduitdeapplgg})}
\label{Sec App 2}

Since $\tilde T_n(x)\geq\e_n$ for all $x\in\bR^d$, we have 
\begin{eqnarray*}
\lefteqn{\limsup_{n\tinf} \frac{v_n^2}{nh_n^{*d}} \log
\bP \, \left[ \, \sup_{x\in C_n^c} 
  \frac{{v_n}\left|f_n^*(x)-f(x)\right|}{\sqrt{\tilde T_n(x)\ka}}
  \, > \,  \delta \, \right]}\\
& \leq &
\limsup_{n\tinf} \frac{v_n^2}{nh_n^{*d}} \log
\bP \, \left[ \, \sup_{\|x\|>w_n} 
  \frac{{v_n}\left|f_n^*(x)-f(x)\right|}{\sqrt{\e_n\ka}}
  \, > \,  \delta \, \right].
\end{eqnarray*}
Now, recall that  
$\sup_{\|x\|\in\bR^d}|\bE(f_n^*(x))-f(x)|=O(h_n^{*2})$; in 
view of the condition $v_nh_n^{*2}/\sqrt{\e_n}\ra 0$ 
in (\ref{assumption}), Equation (\ref{ancien deduitdeapplgg})
is thus a straightforward consequence of the asymptotic
\begin{equation}
\label{applgg}
  \limsup_{n\tinf} \frac{v_n^2}{nh_n^{*d}} 
\log \bP \, \left[ \,  \sup_{\|x\|\geq w_n}
\frac{v_n|f_n^*(x)-\bE(f_n^*(x))|}{\sqrt{\e_n\kappa}}
 \geq \delta  \, \right] \; = \; -\infty.
\end{equation}
To prove (\ref{applgg}), we need the following technical Lemma.

\begin{lemma}
\label{majgg} 
Assume that (A1) and (A4) hold. For all $\gamma>0$, we have
$$
 \sup_{\|x\|\geq w_n}\int_{\bR^d}K^2\left(\frac{x-y}{h_n^*}\right)f(y)dy
 \; \leq \; \gamma h_n^{*d}\e_n.$$
\end{lemma}
\vspace{0.5cm}

\paragraph{Proof of Lemma \ref{majgg}} 
Set  $\gamma>0$, and write
\begin{eqnarray*}
\lefteqn{
\frac{1}{h_n^{*d}\e_n}\int_{\bR^d}K^2\left(\frac{x-y}{h_n^*}\right)f(y)dy}\\
& = & \frac{1}{\e_n}\int_{\|z\|\leq w_n/2}K^2(z)f(x-h_n^*z)dz
+\frac{1}{\e_n}\int_{\|z\|> w_n/2}K^2(z)f(x-h_n^*z)dz
\end{eqnarray*}
On the one hand, we note that
\begin{eqnarray*}
\|x\|\geq w_n\ \ \mbox{and}\ \ \|z\|\leq w_n/2 & \imp &
\|x-h_n^*z\|\geq w_n |1-h_n^*/2|\\
 & \imp &
\|x-h_n^*z\|\geq w_n /2 \ \ \mbox{for $n$ large enough}\\
 & \imp &
f(x-h_n^*z)\leq M_f2^q w_n^{-q}\ \ \mbox{for $n$ large enough}
\end{eqnarray*}
(where $M_f=\sup_{x\in\bR^d} \|x\|^{q}f(x)$), so that
\begin{eqnarray*}
 \sup_{\|x\|\geq w_n}
\frac{1}{\e_n}\int_{\|z\|\leq w_n/2}K^2(z)f(x-h_n^*z)dz
& \leq & 
\frac{M_f2^q}{w_n^q\e_n}\int_{\bR^d}K^2(z)dz\\
& \leq & 
\frac{\gamma}{2}\ \ \mbox{for $n$ large enough}
\end{eqnarray*}
(since $w_n^q\e_n\tinf$).\\
On the other hand, we have
\begin{eqnarray*}
 \sup_{\|x\|> w_n}
\frac{1}{\e_n}\int_{\|z\|> w_n/2}K^2(z)f(x-h_n^*z)dz
& \leq & 
\sup_{\|x\|> w_n}
\frac{\|f\|_\infty 2^q}{w_n^q\e_n}\int_{\bR^d}\|z\|^q K^2(z)dz\\
& \leq & 
\frac{\gamma}{2}\ \ \mbox{for $n$ large enough}
\end{eqnarray*}
which concludes the proof of Lemma \ref{majgg}.
$\Box$\\
\\
Let us now come back to the proof of (\ref{applgg}).
Set 
$${\cal F}_n=\left\{K\left(\frac{x-.}{h_n^*}\right);\ \|x\|> w_n\right\}.$$
The classes ${\cal F}_n$, $n\geq 1$, are contained in the class 
${\cal F}(K)$ defined by (\ref{Fcovering}); since $K$ satisfies the 
covering number condition (\ref{covering}), there exist $A>0$ and $v>0$
such that, $\forall\e>0$, $\forall n\geq 1$,
$$N_2(\e\|K\|_\infty ,\bP ,{\cal F}_n)\leq 
{\left(\frac{A}{\e}\right)}^{v}.$$
Now, let us take
$U=\|K\|_\infty$ and  $\sigma^2=\gamma h_n^{*d}\e_n$ with $\gamma>0$.
Since $h_n^{*d}\e_n\ra 0$ and in view of Lemma \ref{majgg}, we have, 
for $n$ sufficiently large, 
$$\sigma\leq U/2$$ 
and
$$
\sigma^2 \; \geq \; \bE\left[K^2\left(\frac{x-X_1}{h_n^*}\right)\right]
\ \ \mbox{for}\ \  {\|x\|\geq w_n}.$$
Thus, it follows from Theorem 2.1 of Gin\'e and Guillou (2002) that there
exist two constants $M$ and $L$ depending only on $A$ and $v$ such that,
for 
$$t\geq L\left(U\log\frac{U}{\sigma}+\sqrt{n}\sigma
\sqrt{\log\frac{U}{\sigma}}\right),$$
we have
\begin{eqnarray}
\lefteqn{\log \bP \, \left[ \,  \sup_{\|x\|\geq w_n}
\left|\sum_{i=1}^n\left[K\left(\frac{x-X_i}{h_n^*}\right)-
\bE\left(K\left(\frac{x-X_i}{h_n^*}\right)\right)\right]\right|
 \geq t  \, \right]}\nonumber\\  
& \leq &
\log M-\frac{t}{MU}\log\left(1+\frac{tU}{M
[\sqrt{n}\sigma+U\sqrt{\log({U}/{\sigma})}]^2
}\right).
\label{ulti}
\end{eqnarray}
It follows from the conditions ${h_n^{*}}/{{\e_n}}\rightarrow 0$,
$v_n\e_n^{3/2}\rightarrow \infty$ and
${nh_n^{*d}}/[{v_n^2\log(1/h_n^*)}]\rightarrow \infty$ in (\ref{assumption})
that
$\lim_{n\tinf}\sqrt{n}\sigma/[U\sqrt{\log\frac{U}{\sigma}}]=\infty$, and thus,
for $n$ large enough,
\begin{equation}
\label{2conc}
\sqrt{n}\sigma+U\sqrt{\log\frac{U}{\sigma}}\leq 2\sqrt{n}\sigma.
\end{equation}
Now, set $t_n=\delta n h_n^{*d}\sqrt{\e_n\kappa}/v_n$.
It follows from (\ref{2conc}) that, for $n$ large enough,
\begin{eqnarray*}
\frac{U\log\frac{U}{\sigma}+\sqrt{n}\sigma
\sqrt{\log\frac{U}{\sigma}}}{t_n}
& \leq & 
\frac{2\sqrt{n}\sigma
\sqrt{\log\frac{U}{\sigma}}}{t_n}\\
& \leq & 
\sqrt{\frac{4v_n^2\log(U/[\gamma\e_nh_n^{*d}])}{\kappa nh_n^{*d}}} .
\end{eqnarray*}
The conditions ${h_n^{*}}/{{\e_n}}\rightarrow 0$ and
${nh_n^{*d}}/[{v_n^2\log(1/h_n^*)}]\rightarrow \infty$ in (\ref{assumption})
ensure then that, for $n$ large enough, 
\begin{equation}
\label{1conc}
t_n\geq L\left(U\log\frac{U}{\sigma}+\sqrt{n}\sigma
\sqrt{\log\frac{U}{\sigma}}\right).
\end{equation}
Since
\begin{eqnarray*}
\lefteqn{\log \bP \, \left[ \,  \sup_{\|x\|\geq w_n}
\frac{v_n|f_n^*(x)-\bE(f_n^*(x))|}{\sqrt{\e_n\kappa}}
 \geq \delta  \, \right]}\\
& = &
\log \bP \, \left[ \,  \sup_{\|x\|\geq w_n}
\left|\sum_{i=1}^n\left[K\left(\frac{x-X_i}{h_n^*}\right)-
\bE\left(K\left(\frac{x-X_i}{h_n^*}\right)\right)\right]\right|
 \geq t_n  \, \right],
\end{eqnarray*}
it follows from (\ref{ulti}), (\ref{2conc}) and (\ref{1conc}) that, 
for $n$ large enough,
\begin{eqnarray*}
\log \bP \, \left[ \,  \sup_{\|x\|\geq w_n}
\frac{v_n|f_n^*(x)-\bE(f_n^*(x))|}{\sqrt{\e_n\kappa}}
 \geq \delta  \, \right]
& \leq &
\log M-\frac{t_n}{MU}\log\left(1+\frac{t_nU}{4Mn\sigma^2}\right).
\end{eqnarray*}
Noting that
$$\lim_{n\tinf}\frac{t_nU}{4Mn\sigma^2}=\lim_{n\tinf}\frac{\delta U}
{4M\gamma v_n\sqrt{\e_n}}=0,$$
we obtain
\begin{eqnarray*}
 \limsup_{n\tinf} \frac{v_n^2}{nh_n^{*d}} 
\log \bP \, \left[ \,  \sup_{\|x\|\geq w_n}
\frac{v_n|f_n^*(x)-\bE(f_n^*(x))|}{\sqrt{\e_n\kappa}}
 \geq \delta  \, \right] 
& \leq &
-\limsup_{n\tinf} \frac{v_n^2}{nh_n^{*d}} \ \frac{t_n^2}{4M^2n\sigma^2}\\
& \leq &
-\frac{\delta^2\kappa}{4\gamma M^2},
\end{eqnarray*}
which implies (\ref{applgg}) by letting $\gamma\rightarrow 0$.


\end{document}